\documentclass[number,citesort,dvips]{arxbj}
\usepackage{graphics}
\aid{0}
\volume{16}
\issue{4}
\pubyear{2010}
\firstpage{1312}
\lastpage{1342}
\doi{10.3150/09-BEJ243}

\makeatletter
\newcommand{\ud}{\mathrm{d}}

\newtheorem{theorem}{Theorem}
\newtheorem{corollary}{Corollary}
\newtheorem{lemma}{Lemma}
\newtheorem{proposition}{Proposition}
\newtheorem{conjecture}{Conjecture}

\newremark{remark}{Remark}

\def\Exp{{\mathbb{E}}}
\def\Pr{{\mathbb{P}}}
\def\R{\mathbb{R}}
\def\Z{\mathbb{Z}}
\def\F{{\mathcal{F}}}
\def\D{\mathcal{D}}
\def\N{\mathbb{N}}
\def\eps{{\varepsilon}}
\def\1{\mathbf{1}}

\makeatother

\begin{document}
\begin{frontmatter}

\title{Passage-time moments and hybrid zones for the exclusion-voter model}
\runtitle{The exclusion-voter model}

\begin{aug}
\author[a]{\fnms{Iain M.} \snm{MacPhee}\thanksref{a,e1}\ead[label=e1,mark]{i.m.macphee@durham.ac.uk}},
\author[a]{\fnms{Mikhail V.} \snm{Menshikov}\thanksref{a,e2}\ead[label=e2,mark]{mikhail.menshikov@durham.ac.uk}},
\author[b]{\fnms{Stanislav}~\snm{Volkov}\thanksref{b}\ead[label=e3]{s.volkov@bristol.ac.uk}} \and
\author[c]{\fnms{Andrew R.} \snm{Wade}\corref{}\thanksref{c}\ead[label=e4]{andrew.wade@strath.ac.uk}}
\runauthor{MacPhee, Menshikov, Volkov and Wade}
\address[a]{Department of Mathematical Sciences, Durham University,
South Road,  Durham DH1 3LE, UK.
E-mails: \printead*{e1,e2}}
\address[b]{Department of Mathematics, University of Bristol, University Walk,
Bristol BS8 1TW, UK.\\ \printead{e3}}
\address[c]{Department of Mathematics and Statistics, University of Strathclyde,
26 Richmond Street,\\
Glasgow G1 1XH, UK. \printead{e4}}
\end{aug}

% HISTORY:
\received{\smonth{10} \syear{2008}}
\revised{\smonth{7} \syear{2009}}

% ABSTRACT
\begin{abstract}
We study the non-equilibrium
dynamics of
a one-dimensional
interacting particle system
that is a mixture of the voter model and the exclusion process.
With the process started from
a finite perturbation of the ground state
Heaviside configuration consisting of $1$'s to the left of the origin and $0$'s
elsewhere, we study the
relaxation time $\tau$, that is, the first hitting time
of the ground state configuration (up  to translation).
We give conditions
for $\tau$ to be finite and for certain moments of $\tau$
to be finite or infinite,
and prove a result that approaches a conjecture of
Belitsky~\textit{et al.}~(\textit{Bernoulli} \textbf{7} (2001) 119--144).
Ours are the first non-existence-of-moments results for $\tau$
for the mixture model.
Moreover, we give
almost sure asymptotics for the
evolution of the size of the hybrid
(disordered) region.
Most of our results pertain to the discrete-time setting, but
several transfer to continuous-time.
As well as the mixture process, some of our results also
cover   pure exclusion.
We state several significant
open problems.
\end{abstract}

% KEYWORDS

\end{frontmatter}

%s1 ###
\section{Introduction}

The \textit{exclusion-voter} model
is a one-dimensional lattice-based
interacting particle
process with nearest-neighbour interactions, introduced by
Belitsky \textit{et al.}~in \cite{bfmp}, that is, a mixture of the
symmetric voter model and the simple exclusion process. For
background on the latter two models (separately)
and interacting particle systems in general, see
\cite{ligg,ligg2}.

The voter model has been used to
model  the spread of an
opinion through a static population via
nearest-neighbour interactions; see, for example,~\cite{hl}.
The mixture model studied here is a
natural extension of this model whereby individuals
do not have to remain static, but may move by
switching places. Alternative motivations, such as from the point
of view of competition of species (see, e.g.,~\cite{cliff}) can also
be adapted to the mixture model.
As our results show, allowing place-swaps
can have a dramatic effect on the dynamics of the process.

The exclusion-voter model is a Markov  process with
state space $\{0,1\}^\Z$;
each site of   $\Z$
can be labelled
either $0$ or $1$,  representing the presence of one of two types of particle.
The ground state of our model will be the `Heaviside' configuration
$\ldots 111 000\ldots.$
We consider initial configurations that are finite
perturbations of this ground state and so
contain
a finite number of unlike pairs, where, by `pair', we always mean two adjacent particles.

In this paper, we concentrate on a \textit{discrete-time} process
that can be described informally as follows.
At each time step, the \textit{simple exclusion process} selects uniformly at random from amongst
all unlike pairs. If the chosen pair is $01$, it flips  to $10$ with probability $p$ (otherwise there is no change);
if the pair is
$10$, it flips to $01$ with probability $1-p$. On the other hand,
at each time step, the \textit{symmetric voter model} selects uniformly at
random from all unlike pairs and then flips the chosen pair to
either
$00$ or $11$, with equal chance of each. The model that is
considered in this paper, introduced in \cite{bfmp}, is a mixture of
these two processes where, at each time step, we determine
independently at random
whether to perform a voter-type move (with probability $\beta$)
or an exclusion-type move (probability $1-\beta$).

The analogous continuous-time exclusion-voter
model can be defined via its infinitesimal
generator and
constructed via a Harris-type
graphical construction. The
discrete-time process described above
is naturally embedded in the continuous-time
process. In our analysis, we work in discrete time,
and the discrete-time process has its own interest,
but, as we shall indicate,
some of our results
transfer almost immediately into continuous time.

Individually, the exclusion process and voter model exhibit very different
behaviour. For instance, in the exclusion process, there is local
conservation of $1$'s: the number of $1$'s
in a bounded interval can change only through the boundary. There is no such
conservation in the voter model.
In the mixture process that we study in the present paper,
voter moves and exclusion moves interact in a highly non-trivial way. This introduces
technical difficulties: for instance, voter moves can
cause drastic changes quickly, also there is no obvious monotonicity
property.
We describe the model more formally and state our results in the next section. First, we
outline the existing literature and
the contribution of the present paper.

In \cite{bfmp},  results were proven for the exclusion process and voter model separately,
as well as some initial results for the mixture model. The main problems
left open in \cite{bfmp} were the non-existence of passage-time moments
and the issue of transience/recurrence for the mixture model.
As we will describe shortly, the present paper makes contributions
to each of these problems.
Some of the results in \cite{bfmp}, in the symmetric exclusion ($p=1/2$)
case,
are generalized to non-nearest-neighbour
interactions in \cite{ss}. Certain `ergodic' properties
of a generalization of the continuous-time exclusion-voter model, again in the symmetric exclusion case,
are studied in \cite{jung}. The goal of the present paper is
to study the mixture model in  more depth than \cite{bfmp}. In particular, we prove
new results on: (i) the passage-time problem for the exclusion-voter
model, the main contribution being the
(more difficult) non-existence of passage-time
moments;  (ii) the size of the
disordered region where $1$'s and $0$'s intermingle.
This
region we call the \textit{hybrid zone} (cf.~\cite{coxdurrett}).
Our results  leave several  open
problems and we put forward some conjectures with
regard to these in the next section.

Let us describe more specifically
the contribution of the present
paper to the passage-time problem
for the exclusion-voter model.
The passage time of interest to us here is
the \textit{relaxation time} $\tau$
-- the return time of the configuration to the ground state.
In general,
one can
often prove the existence of  moments of passage times
directly via semimartingale
(Lyapunov-type function) criteria such as those in
\cite{lamp2,aim,ai2}, in the vein of Foster \cite{foster}.
The
non-existence of moments (for which no results
have previously been obtained for the exclusion-voter model with $\beta \in (0,1)$)
is usually a harder problem.
In general,
semimartingale-type arguments are available in this case too (see, e.g.,~\cite{lamp2,aim,ai1}), but under more restrictive conditions
than the corresponding existence results: non-existence results
typically need fine control over jumps of the process. Lamperti \cite{lamp2}
was first to establish a general methodology for proving
non-existence of passage-time moments, based on
finding a suitable submartingale and obtaining a good-probability
lower bound for passage times; his
method was later extended in \cite{aim,ai1}.
The same two elements
form the basis of our approach, but we must proceed
differently since the exclusion-voter model
does not possess the regularity required
by existing general results such as those
of \cite{lamp2,aim,ai1}.

On the one hand, we extend the region of the parameter space
of the model for which almost sure finiteness of $\tau$ is known and we give
results on the existence of higher
moments of $\tau$
(including in the case of pure exclusion). On the other hand, we show
the \textit{non-existence}
of certain moments of $\tau$;
this problem was not addressed in \cite{bfmp}.
Each of these opposing directions requires
us to develop new techniques.
We prove, for example, that under
certain conditions, $1+\eps$ moments ($\eps>0$)
of $\tau$ do not exist; this
approaches a  conjecture in \cite{bfmp}.

The second main contribution of the paper is
to study the evolution of the size of the hybrid zone. Our
basic tools are again  semimartingales: we apply  general results
on   almost sure bounds for stochastic processes
from \cite{mvw}.
For instance, for the pure exclusion
process in the case $p=1/2$ we prove
that, with probability 1,
the maximum size
of the hybrid zone
up to time $t$ remains
bounded between $t^{1/3}$ and $t^{1/2}$,
ignoring logarithmic factors.

In the next section, we give
some more formal definitions, state our main
results and discuss some (challenging)
open problems.

%s2 ###
\section{Definitions and statement of results}\label{results}

We now formally
describe the model that we study, as considered in \cite{bfmp}.
We introduce some notation to describe the configuration of the process.
Let $\D' \subset \{0,1\}^\Z$ denote the set of configurations with a finite
number of $0$'s to the left of the origin and $1$'s to the right. Let `$\sim$' denote the equivalence
relation on $\D'$ such that for $S, S' \in \D'$, $S \sim S'$ if and only
if $S$ and $S'$ are translates of each other. Then set $\D:= \D' / \sim$. In other words,
the
configuration space $\D$ is the set of configurations
taking the   form of an infinite string of $1$'s followed by a finite number of $0$'s and $1$'s followed by an infinite string of $0$'s, modulo translations.
For example, one configuration $S \in \D$ is
\begin{equation}\label{configex}
S= \ldots 1110000000011100001001001000000001111000 \ldots.
\end{equation}
Configurations such as those in $\D$ are sometimes called \textit{shock profiles} (see, e.g.,~\cite{bfmp}).

Fix
$\beta \in [0,1]$ (the mixing parameter) and $p \in [0,1]$
(the exclusion parameter).
The discrete-time exclusion-voter process $\xi=(\xi_t )_{ t \in \Z^+}$
with parameters $(\beta, p)$
is a time-homogeneous Markov chain on the countable state space $\D$. The one-step
transition probabilities are determined by the following mechanism.
At each time step, we decide independently at random whether to perform a
\textit{voter} move or an \textit{exclusion} move. We choose a voter move with probability $\beta$ and an exclusion move
with probability $1-\beta$. Having decided this, we choose an unlike adjacent pair (i.e.,~$01$ or $10$)
uniformly at random.
The voter move is such that the chosen pair ($01$ or $10$) flips to $00$ or $11$, each with probability $1/2$.
The exclusion move is such that a chosen pair $01$ flips to $10$ with
probability $p$ (otherwise no move) and a chosen pair $10$ flips to $01$ with probability $q:=1-p$
(otherwise no move).

In addition to the discrete-time model
that is the focus of the present paper, there is a corresponding continuous-time model,
also introduced in \cite{bfmp}.
A priori, the
relationship between the two time-scales is complicated,
but from our results on the discrete-time process,
we can obtain some results in the continuous-time
setting too. For a description of the continuous-time model,
its relationship to the discrete-time model
that is our main object of study and our results in the continuous-time setting,
see Section~\ref{cont} below.

We denote the underlying probability space for $\xi$
by $(\Omega,\F,\Pr_{\beta,p})$ and the corresponding expectation by
$\Exp_{\beta,p}$.
We denote the ground state Heaviside configuration $\D_0 \in \D$, which
consists of a single  pair $10$ abutted by infinite strings of $1$'s and $0$'s to the
left and right, respectively:
\[ \D_0 =  \ldots 11110000 \ldots,\]
up to translation.
The next result gives some
elementary properties of the state space
$\D$ under $\Pr_{\beta,p}$. In particular, Proposition~\ref{prop1} says that for $(\beta,p) \in (0,1)^2$
(i.e., in the interior of the parameter space), $\xi$ is irreducible and aperiodic
under $\Pr_{\beta,p}$.
%p1
\begin{proposition}\label{prop1}
$\D_0$ is an absorbing state under $\Pr_{\beta,1}$ for any $\beta \in [0,1]$.
Suppose that $\beta \neq 1$ and $(\beta,p) \notin \{ (0,0), (0,1) \}$. All states
in $\D \setminus \{ \D_0 \}$ then communicate under $\Pr_{\beta,p}$.
Suppose that  $\beta \neq 1$, $p<1$ and $(\beta,p) \neq (0,0)$. All states
in $\D$ then communicate under $\Pr_{\beta,p}$, and $\xi$ is irreducible and
aperiodic.
\end{proposition}

For $S_0 \in \D$, define the \textit{relaxation time}
for the process $\xi$  as
\[
\tau:= \min \{ t \in \N\dvt  \xi_t = \D_0   \}.
\]
We introduce some convenient terminology.
If
$\Pr_{\beta,p} (\tau = +\infty  |  \xi_0 = S_0 ) > 0$ for
$S_0 \in \D \setminus\{ \D_0\}$, then we say that $\xi$
is \textit{transient} started from $S_0$;
if $\Pr_{\beta,p} (\tau <\infty  |  \xi_0 = S_0 ) =1$ for
$S_0 \in \D \setminus\{ \D_0\}$, then we say that $\xi$
is \textit{recurrent} started from $S_0$. In the latter
case, if, in addition,
$\Exp_{\beta,p} [ \tau   |  \xi_0 =S_0 ] < \infty$
for
$S_0 \in \D \setminus\{ \D_0\}$,
then we say that $\xi$ is \textit{positive recurrent} started from $S_0$.
When $\xi$ is irreducible  (see Proposition~\ref{prop1}),
this terminology coincides with the standard usage
for countable state space Markov chains. When $\xi$ is irreducible and aperiodic
(see Proposition~\ref{prop1}),
we may use the term
\textit{ergodic} in the positive recurrent case.

Results of Liggett (see, e.g.,~Chapter VIII of \cite{ligg}) imply that  the pure exclusion process ($\beta=0$)
is positive recurrent for all $S_0 \in \D$
if and only if $p>1/2$.
We recall
the following result, which
is contained in Theorems 5.1, 5.2, 6.1, 7.1 and 7.2 of \cite{bfmp}, together
with an inspection of (7.2) in \cite{bfmp}
for part (iii)(a).
%t1
\begin{theorem}\label{bfmpthm}
\textup{(i)} Suppose that $\beta=0$ (pure exclusion). Then, for
any $S_0 \in \D$, $\xi$ is positive recurrent
for $p>1/2$ and transient for $p \leq 1/2$.
\smallskipamount=0pt
\begin{longlist}[(iii)]
\item[(ii)] Suppose that $\beta=1$ (pure voter). Then $\xi$ is positive recurrent
for any $S_0 \in \D$
and, moreover, for any $S_0 \in \D \setminus\{ \D_0 \}$ and any $\eps>0$,
\[
\Exp_{1,p} \bigl[ \tau ^{(3/2)-\eps}  |  \xi_0 =S_0 \bigr] < \infty;\qquad
\Exp_{1,p} \bigl[ \tau ^{(3/2)+\eps}  |  \xi_0 = S_0 \bigr]
=\infty.
\]
\item[(iii)] Suppose that $\beta \in (0,1)$ (mixture process).
\begin{enumerate}[(iii) (a)]
\item[(a)]
If $\beta$ and $p \in [0,1]$ are such that $(1-p)(1-\beta) < 1/3$,
then $\xi$ is positive recurrent for any $S_0 \in \D$.
In particular, for any $\beta > 2/3$ and any $p \in [0,1]$,
$\xi$ is positive recurrent for any $S_0 \in \D$.
\item[(b)]
For $p \geq 1/2$ and any $\beta>0$, $\xi$ is positive recurrent for any $S_0 \in \D$.
\end{enumerate}
\end{longlist}
\end{theorem}

In \cite{bfmp}, the following was Conjecture 7.1.
%c1
\begin{conjecture}\label{conj1}
For any $p < 1/2$, there exists $\beta_0=\beta_0(p)>0$ such that for any $\beta < \beta_0$,
$\xi$ is not positive recurrent,
that is,~$\Exp_{\beta,p} [ \tau  |  \xi_0 = S_0 ] = \infty$ for any $S_0 \in \D \setminus\{ \D_0 \}$.
\end{conjecture}

Our first result says that for $\beta$ small enough (so that the
exclusion part is prevalent), $1+\eps$ moments do not exist;
thus Conjecture~\ref{conj1} remains tantalizingly open.
%t2
\begin{theorem}\label{thm1}
For each $p < 1/2$, there exists   $\beta_1 = \beta_1 (p) = (1-2p)/(2-2p) \in (0,1/2]$ such that for
all $\beta \leq \beta_1$,  any $\eps>0$ and any $S_0 \in \D \setminus\{ \D_0 \}$,
\[ \Exp_{\beta,p} [ \tau^{1+\eps}   |  \xi_0 = S_0 ] = \infty.\]
\end{theorem}

Our second result
says that in the mixture process, the
presence of a  transient exclusion ensures that $2+\eps$ moments do not exist.
Thus, for $p \leq 1/2$, even in the case where Theorem~\ref{bfmpthm}(iii)
applies,
the  recurrence is polynomial in nature, that is,~`heavy tailed'.
%t3
\begin{theorem}\label{thm2} Suppose that $p \leq 1/2$, $\beta \in [0,1]$.
For any $\eps>0$ and  $S_0 \in \D \setminus\{ \D_0 \}$,
\[ \Exp_{\beta,p} [ \tau ^{2+\eps}  |  \xi_0 = S_0 ] = \infty.\]
\end{theorem}

In view of Theorem~\ref{bfmpthm}(ii), we
suspect that mixing transient $(p\leq 1/2)$ exclusion
with the voter model ought not to lead to a lighter tail for $\tau$, as we now conjecture.
%c2
\begin{conjecture}\label{conj4} Suppose that  $p \leq 1/2$, $\beta \in [0,1]$. For any $\eps>0$ and  $S_0 \in \D \setminus\{ \D_0 \}$,
\[ \Exp_{\beta,p} \bigl[ \tau ^{(3/2)+\eps}  |  \xi_0 = S_0 \bigr] = \infty.\]
\end{conjecture}

Even this conjecture seems to be challenging, as   exclusion
and voter moves interact in complex ways. Technically, the issue
that prevents us from reducing the $2$ to $3/2$ in Theorem~\ref{thm2}
is that  exclusion moves can (and typically will)
increase the number of blocks.

An open problem mentioned in \cite{bfmp}
is whether the mixture process with
$\beta>0$ and
$p<1/2$  is, in fact, transient (it is recurrent
for $p \geq 1/2$, by Theorem~\ref{bfmpthm}(iii)(b)).
Simulations that we have performed have been inconclusive.
We conjecture the following.
%c3
\begin{conjecture}\label{conj2}
Suppose that $p<1/2$, $\beta >0$. For any $S_0 \in \D$,
$\xi$ is recurrent.
\end{conjecture}

Note that if Conjectures~\ref{conj1} and~\ref{conj2}
both hold, then there is \textit{null recurrence}
for $p < 1/2$ and $\beta \in (0, \beta_0)$.
Our next result represents some progress
in the direction of Conjecture~\ref{conj2} and gives
recurrence in a previously unexplored
region
of the parameter space.
%t4
\begin{theorem}\label{recthm}
Suppose that
$p < 1/2$, $\beta \geq 4/7$.
For any $S_0 \in \D$,
$\xi$ is recurrent.
\end{theorem}

We now turn to the problem of existence of moments for $\tau$. First, we
consider the pure exclusion process in the positive recurrent ($p>1/2$) case.
If we further restrict to $p>2/3$, then it is possible to construct a positive
strict supermartingale with uniformly bounded increments (see (5.7) in \cite{bfmp})
and so it is not hard to show that
all polynomial moments of $\tau$ exist in that case. Theorem~\ref{thmex1} below extends
this conclusion to all $p>1/2$.
%t5
\begin{theorem}\label{thmex1}
Suppose that $\beta =0$, $p>1/2$. For any $S_0 \in \D$
and any $s \in [0,\infty)$,
\[ \Exp_{0,p} [ \tau^{s}  |  \xi_0 = S_0 ] < \infty.\]
\end{theorem}

We suspect that under the conditions of
Theorem~\ref{thmex1}, the existence
of some superpolynomial `moments' for $\tau$ can be
obtained via our techniques and general results from \cite{ai1}.
The next result covers the mixture process in the case where the
exclusion component is positive recurrent.
In the $\beta \in [0,1]$, $p>1/2$ case, we know from Theorem~\ref{bfmpthm}
that $\Exp [ \tau ]<\infty$; the next theorem says that some higher moments
are also finite.
%t6
\begin{theorem}\label{thmex2}
Suppose that $\beta \in [0,1]$, $p>1/2$. For any $S_0 \in \D$,
$\Exp_{\beta,p} [ \tau^{6/5}  |  \xi_0 = S_0 ] < \infty$.
\end{theorem}

In view of Theorem~\ref{bfmpthm} and Theorem~\ref{thmex1}, in the setting
of Theorem~\ref{thmex2}, we are mixing together the voter model, for which $(3/2)-\eps$
moments exist, and the recurrent exclusion process, for which all moments exist.
One might therefore  hope to improve the exponent in Theorem~\ref{thmex2}
to at least $(3/2)-\eps$;  this is another open problem.

\begin{figure}

\includegraphics{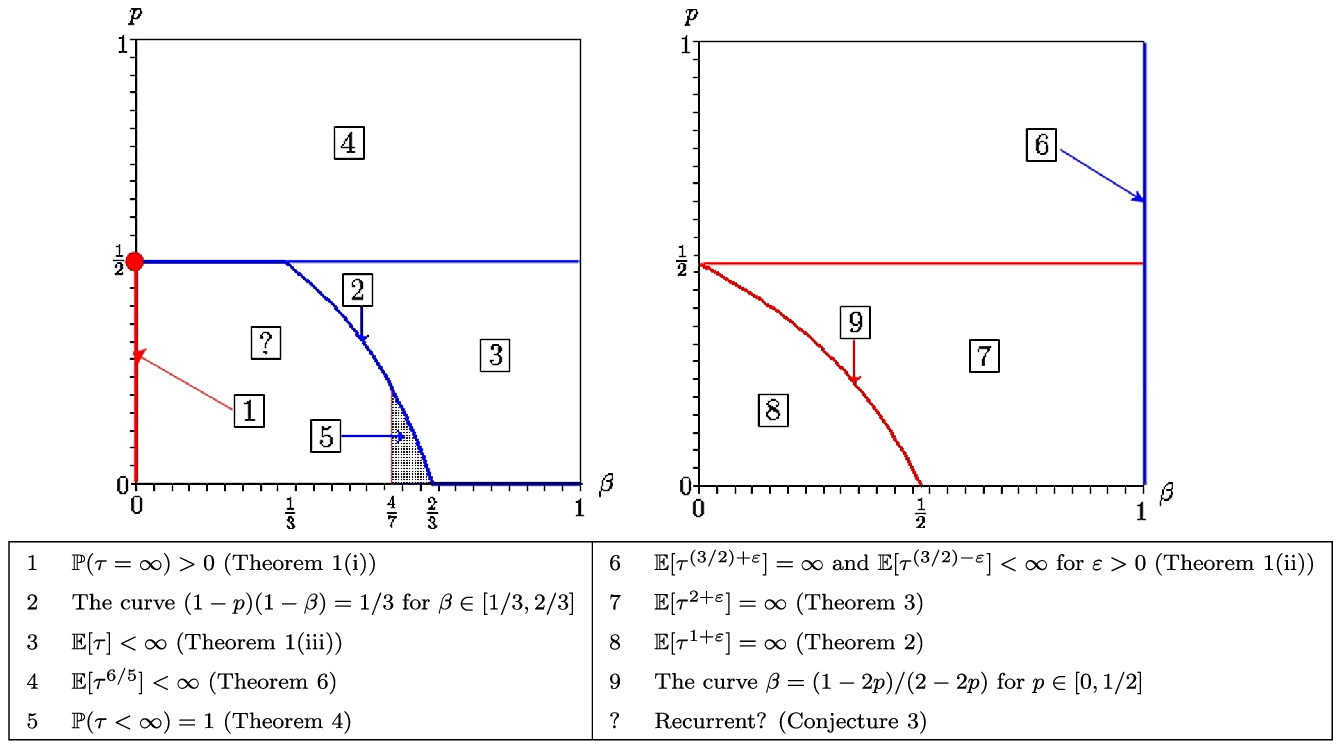}

%1 &  $\Pr(\tau=\infty)>0$ (Theorem~\ref{bfmpthm}(i)) & 6 &
% $\Exp[\tau^{(3/2)+\eps}]=\infty$ and $\Exp[ \tau^{(3/2)-\eps}]< \infty$
%for $\eps>0$ (Theorem~\ref{bfmpthm}(ii))  \\
%2 & The curve $(1-p)(1-\beta)=1/3$ for $\beta \in [1/3,2/3]$ & 7 & $\Exp [\tau^{2+\eps}]=\infty$ (Theorem~\ref{thm2}) \\
%3 & $\Exp[\tau]<\infty$ (Theorem~\ref{bfmpthm}(iii)) & 8 & $\Exp[ \tau^{1+\eps}] =\infty$ (Theorem~\ref{thm1}) \\
%4 & $\Exp[\tau^{6/5}]<\infty$ (Theorem~\ref{thmex2}) & 9 &The curve $\beta= (1-2p)/(2-2p)$ for $p \in [0,1/2]$\\
%5 & $\Pr(\tau < \infty)=1$ (Theorem~\ref{recthm}) & ? &Recurrent? (Conjecture~\ref{conj2})  \\
%}
\caption{Representations of the $(\beta,p)$ parameter space. The   key  explains the labelling, together
with the appropriate result from the text (for brevity, we have dropped the subscripts from
$\Pr,\Exp$).
}\label{fig1}
\end{figure}

Figure~\ref{fig1} gives two diagrams
of the $(\beta,p)$ parameter space,
summarizing the results of the previous
theorems for the relaxation time $\tau$.

%f1 ###

We now state our results on the size of the hybrid zone.
First, we need to
introduce some more notation, following \cite{bfmp}.
A $1$-block ($0$-block)
is a maximal string of consecutive $1$'s ($0$'s). Configurations in $\D$
consist of a finite number of such blocks.
For
$S \in \D$, let $N=N(S)\ge 0$ denote the number of $1$-blocks not
including the  infinite $1$-block to the left (this is the same as
number of $0$-blocks, not including the   infinite
$0$-block to the right). Enumerating left-to-right,
let $n_i = n_i (S)$ denote the size of the $i$th
$0$-block and $m_i=m_i(S)$ the size of the $i$th $1$-block. We may represent configuration
$S \in \D \setminus \{ \D_0\}$ by the vector
$(n_1,m_1,\ldots, n_N,m_N)$. For example, the configuration $S$
of (\ref{configex}), which has $N(S)= 5$,
has the representation $(8,3,4,1,2,1,2,1,8,4)$.
Set $|\D_0|:=0$ and, for $S \in \D \setminus \{ \D_0\}$, let
$|S|:=\sum_{i=1}^N (n_i+ m_i)$ represent the size
of the hybrid zone, that is, the length of the string of $0$'s and $1$'s between the infinite string
of $1$'s to the left and the infinite string of $0$'s to the right.

The next result gives upper bounds
for the size of the hybrid zone $|\xi_t|$
and the number of blocks $N(\xi_t)$;
in particular,
part (ii) covers the case $\beta=0$, $p=1/2$ of the symmetric pure (transient) exclusion process.
%t7
\begin{theorem}\label{thm3}
\textup{(i)}
Suppose that $\beta \in [0,1]$, $p \in [0,1]$. For any $\eps>0$, $\Pr_{\beta,p}$-a.s.,
for all but finitely many $t \in \Z^+$,
\begin{equation}\label{yy0}
\max_{0 \leq s \leq t} N(\xi_s ) \leq
\cases{
t^{1/2} (\log t)^{(1/2)+\eps}, &\quad\mbox{if }$p < 1/2$, \cr
t^{1/3} (\log t)^{(1/3)+\eps}, &\quad\mbox{if }$p \geq 1/2$.
}
\end{equation}
\smallskipamount=0pt
\begin{longlist}[(ii)]
\item[(ii)]
Suppose that $\beta \in [0,1]$, $p \geq 1/2$. For any $\eps>0$, $\Pr_{\beta,p}$-a.s.,
for all but finitely many $t \in \Z^+$,
\begin{equation}\label{yy1}
\max_{0 \leq s \leq t} |\xi_s| \leq t^{1/2} (\log t)^{(1/2)+\eps}.
\end{equation}
\end{longlist}
\end{theorem}

The remainder of our results deal with the pure exclusion process
($\beta=0$). In the continuous-time setting,
related results on
the growth of the hybrid zone of the pure exclusion process were
first obtained by Rost \cite{rost} in the totally asymmetric case;
see Section VIII.5 of \cite{ligg}, and
\cite{afs} for more general results. In particular,
Theorems 5.2, 5.3, and 5.12 on pages 403--407 of \cite{ligg} say, very loosely,
that under $\Pr_{0,p}$,
\[
| \eta_t | \approx t \qquad (p < 1/2); \qquad | \eta_t | \approx t^{1/2} \qquad (p = 1/2),
\]
where $\eta$ is the continuous-time version of $\xi$, as
described in Section~\ref{cont}.
In particular, the symmetric case is significantly different
from the asymmetric case.
However,
there seems to be no immediate
way to translate these
results between the continuous- and discrete-time
settings (see Section~\ref{cont} below).
Part (i) of the next result
strengthens the bound in (\ref{yy0})
slightly in the pure exclusion case with $p< 1/2$. Part (ii)
complements the $\beta =0$ case of (\ref{yy1})
for the case $p < 1/2$ (transient but not symmetric exclusion);
it quantifies the rate of transience.
%t8
\begin{theorem}\label{twothirds} Suppose that $\beta =0$ and $p \in
[0,1]$.
\begin{longlist}[(ii)]
\item[(i)]
There exists $C \in (0,\infty)$ such that
for any $p \in [0,1]$, $\Pr_{0,p}$-a.s.,
for all  $t \in \Z^+$,
\[ \max_{0 \leq s \leq t} N(\xi_s) \leq C t^{1/2}.\]
\item[(ii)]
Suppose that $p \in [0,1/2)$. Then,
for any $\eps>0$,
$\Pr_{0,p}$-a.s., for all
but finitely many $t\in\Z^+$,
\[ \max_{0 \leq s \leq t} |\xi_s| \leq t^{2/3} (\log t)^{(1/3)+\eps}.\]
On the other hand, there exists $c(p) \in (0,\infty)$ such that
for any $c \in (0,c(p))$,
$\Pr_{0,p}$-a.s., for all but finitely many $t\in\Z^+$,
\[  |\xi_t| \geq c t^{1/2}.\]
\end{longlist}
\end{theorem}

Our next result complements (\ref{yy1}) in the case
$\beta=0$, $p=1/2$.
%t9
\begin{theorem}\label{thm5}
Suppose that $\beta=0$, $p=1/2$. For any $\eps>0$, $\Pr_{0,1/2}$-a.s.,
for all but finitely many $t \in \Z^+$,
\[ t^{1/3} (\log t)^{-(1/3)-\eps} \leq \max_{0 \leq s \leq t} |\xi_s| \leq t^{1/2} (\log t)^{(1/2)+\eps}.\]
\end{theorem}

It is an open problem to obtain sharper versions of the above results on $|\xi_t|$.
In the pure exclusion $(\beta=0)$ case, we conjecture the following.
%c4
\begin{conjecture}
Suppose that $\beta=0$. If $p<1/2$, then for any $\eps>0$,
$\Pr_{0,p}$-a.s.,
for all but finitely many $t\in\Z^+$, $| \xi_t | \leq t^{(1/2)+\eps}$.
If $p=1/2$, then for any $\eps>0$,
$\Pr_{0,1/2}$-a.s.,
for all but finitely many $t\in\Z^+$,
$t^{(1/3)-\eps} \leq | \xi_t | \leq t^{(1/3)+\eps}$.
\end{conjecture}

The structure of the remainder of the paper is as follows.
In Section~\ref{cont},
we describe the \textit{continuous-time}
version of the exclusion-voter model, how it relates
to the discrete-time version
studied here and which
results can be transferred
without too much extra work.
Section~\ref{prep} contains
preliminary results. In
Section~\ref{speeds},
we collect
general semimartingale
results that we apply in the paper. In Section~\ref{config}, we introduce notation and a
convenient representation
for configurations of the model, and we prove Proposition~\ref{prop1}. In Section~\ref{lemmas}, we give some lemmas
on the Lyapunov-type functions
that we will use throughout the paper.
In Section~\ref{moments},
we prove Theorems~\ref{thm1} and~\ref{thm2}
on passage-time
moments, via a series of lemmas.
In Section~\ref{recprf},
we
prove Theorem~\ref{recthm}. In Section~\ref{secex}, we prove
Theorems~\ref{thmex1}
and~\ref{thmex2}.
In Section~\ref{size}, we
prove Theorems~\ref{thm3},~\ref{twothirds} and~\ref{thm5} on the
size of the hybrid zone and number of blocks.

%s3 ###
\section{Continuous time}\label{cont}

The
exclusion-voter model
may also be defined and studied in continuous time.  First, we recall the definition, following \cite{bfmp}.
Let $\nu = ( \nu(x))_{x \in \Z} \in \{0,1\}^\Z$
so that $\nu(x)$ is the label
($0$ or $1$) at $x$.
For $x, y, z \in \Z$, denote
\begin{eqnarray*}
\nu_{x,y} (z):=
\cases{
\nu(y), &\quad\mbox{if }$z=x$,\cr
\nu(x), &\quad\mbox{if }$z=y$, \cr
\nu(z), &\quad\mbox{if }$z \neq x, y$;
}\qquad
\nu_x (z):=
\cases{
1-\nu(z), &\quad\mbox{if }$z=x$, \cr
\nu(z),   &\quad\mbox{if }$z \neq x$.
}
\end{eqnarray*}
In words, $\nu_{x,y}$ is $\nu$ with labels at $x, y$
interchanged and $\nu_x$ is $\nu$ with the label
at $x$ flipped (i.e., replaced by its opposite).
We introduce Markovian generators $\Omega^{\mathrm{e}}_p$ $(p \in [0,1])$
and $\Omega^{\mathrm{v}}$, defined by their action on functions $f$
on $\{0,1\}^\Z$, by
\begin{eqnarray*}
\Omega^{\mathrm{e}}_p f(\nu) & = &\sum_{x,y } p(x,y) \nu (x) \bigl(1- \nu (y)\bigr) [ f(\nu_{x,y})- f(\nu) ]\quad   \mbox{and} \\
\Omega^{\mathrm{v}} f(\nu) & = &\sum_{x} c(x,\nu)  [ f(\nu_{x})- f(\nu) ],
\end{eqnarray*}
where
$p(x,x-1) =p$, $p(x,x+1)=1-p$ and $p(x,y)=0$ for $| x-y | \neq 1$, and
\begin{eqnarray*}
c(x, \nu):=
\cases{\frac{1}{2}  \bigl( \nu (x-1) +\nu(x+1)  \bigr),
&\quad\mbox{if }$\nu(x) = 0$, \vspace*{2pt}\cr
\frac{1}{2}  \bigl( 2 - \nu (x-1) -\nu(x+1)  \bigr),
&\quad\mbox{if }$\nu(x) = 1$.
}
\end{eqnarray*}

The continuous-time exclusion-voter model
with mixing parameter $\beta \in [0,1]$
and exclusion parameter $p \in [0,1]$
is a Markov process $(\eta'_t)_{t \geq 0}$
on $\mathcal{D}' \subset \{ 0,1\}^\Z$ with generator
$(1-\beta) \Omega^{\mathrm{e}}_p + \beta \Omega^{\mathrm{v}}$.
This induces a Markov process $\eta = (\eta_t)_{t \geq 0 }$
on the space of equivalence classes $\mathcal{D}$ by
taking $\eta_t$ to be the $\sim$-equivalence class
of $\eta'_t$.
The process $\eta$ can be constructed
from an array of homogeneous one-dimensional
Poisson processes via a Harris-type graphical
construction; see page~9 of \cite{bfmp} for details.
With the definitions
in Section~\ref{results} and this section, $\xi$ may be
embedded in $\eta$ in the standard way; again, see \cite{bfmp}.

In the continuous-time setting, the relaxation time is
\[
\tau_{\mathrm{c}}:= \inf \{ t \geq 0\dvt  \eta_t = \D_0 \}.
\]
The natural question is:
given the results in Section~\ref{results}
on $\tau$, what is it possible to say about $\tau_{\mathrm{c}}$?
We now outline which of our discrete-time results
for $\xi$ can be readily
transferred to continuous-time  results
for $\eta$ (cf.~Section~8 of \cite{bfmp}).

First,
as pointed out in \cite{bfmp},
recurrence and transience transfer directly:
\[
\Pr_{\beta,p} ( \tau < \infty  |  \xi_0 = S_0 ) = 1
\quad\iff\quad \Pr_{\beta,p} ( \tau_{\mathrm{c}} < \infty   |  \eta_0 = S_0)=1.
\]
To draw conclusions about   moments
(i.e., tails) of the relaxation times,
it is necessary to know about the comparative  rates
of the two processes.
The
transition rate of the continuous-time
process is, roughly speaking,  proportional to
the
number of blocks so the continuous-time
process tends to evolve at least as fast
as the discrete-time
process.

The pure voter model ($\beta=1$)
is well behaved, in the sense that
it cannot increase the number of blocks. Thus, roughly speaking,
the discrete and continuous
time-scales are directly comparable and results
are more easily transferred. This intuition is formalized
in Section 8 of \cite{bfmp}, where it is shown that for any $s>0$,
\[
\Exp_{1,p} [ \tau^s |  \xi_0 = S_0 ] < \infty \quad\iff\quad \Exp_{1,p} [ \tau_{\mathrm{c}} ^s  |  \eta_0 = S_0]< \infty.
\]
In the general case, without more information on the number of blocks, only one-sided results
are possible a priori. It is shown in Section 8 of \cite{bfmp} that for any $s>0,$
\[
\Exp_{\beta,p} [ \tau^s |  \xi_0 = S_0 ] < \infty \quad\Longrightarrow\quad
\Exp_{\beta,p} [ \tau_{\mathrm{c}} ^s  |  \eta_0 = S_0]< \infty.
\]

Theorem~\ref{bfmpthm} above (proved in \cite{bfmp}) therefore
transfers directly to continuous time and holds
with $\tau_{\mathrm{c}}$ instead of $\tau$;
this is Theorem 1.1 in \cite{bfmp}.
In particular, the $\beta=1$ case of this result
shows that the continuous-time
pure voter model is positive recurrent, a result that
goes back to Cox and Durrett (Theorem 4 of \cite{coxdurrett});
for further study of voter model interfaces and
some generalizations, see \cite{bmsv,bmv,dz}.
Moreover our Theorems~\ref{recthm},~\ref{thmex1}
and~\ref{thmex2} also carry across and
hold with $\tau_{\mathrm{c}}$, yielding the following
corollary.
%cor1
\begin{corollary}\label{recthm2} For any $S_0 \in \D$, we have the following:
\begin{longlist}[(iii)]
\item[(i)]
if
$p < 1/2$, $\beta \geq 4/7$, then $\eta$ is recurrent, that is, $\Pr_{\beta,p} (\tau_{\mathrm{c}} <\infty  |  \eta_0 = S_0) = 1$;
\item[(ii)] if   $\beta =0$, $p>1/2$, then for any $s \in [0,\infty)$, $\Exp_{0,p} [ \tau_{\mathrm{c}}^{s}  |  \xi_0 = S_0 ] < \infty$;
\item[(iii)]
if $\beta \in [0,1]$, $p>1/2$, then $\Exp_{\beta,p} [ \tau_{\mathrm{c}}^{6/5}  |  \xi_0 = S_0 ] < \infty$.
\end{longlist}
\end{corollary}

Corollary~\ref{recthm2}(ii) says that for the standard
(continuous-time) recurrent exclusion process, all moments
of $\tau_{\mathrm{c}}$ exist. This fact may  be known, but we
could not find a reference.

%s4 ###
\section{Preliminaries}\label{prep}

%s4.1 ###
\subsection{Technical tools}\label{speeds}

In this section,
we state some
general  martingale-type results
that we will need. In particular,
we will recall some criteria for obtaining
upper and
lower almost sure bounds for discrete-time
stochastic processes on the half-line given in \cite{mvw}.

Let $(\mathcal{F}_t)_{t \in \Z^+}$ be a filtration
on a probability space $(\Omega, \mathcal{F}, \Pr)$. Let ${X}=(X_t)_{t \in \Z^+}$ be a discrete-time
$(\mathcal{F}_t)$-adapted
stochastic process taking values in $[0,\infty)$.
Suppose that $\Pr ( X_0 = x_0 ) =1$ for some $x_0 \in [0,\infty)$.
For the applications in the present paper,
we will, for instance, take $X_t = | \xi_t |$.
The following result combines a maximal inequality (Lemma
3.1 in \cite{mvw}) with an almost sure upper bound
(contained in Theorem 3.2 of \cite{mvw}).
%l1
\begin{lemma}\label{lem13} Let  $B \in (0, \infty)$ be such that, for all $t \in \Z^+$,
\begin{equation} \label{oww3}
\Exp [ X_{t+1}- X_t  |  \F_t ] \leq B\qquad \mbox{a.s.}
\end{equation}
Then:
\begin{longlist}[(ii)]
\item[(i)] for any $r >0$ and any $t \in \N$,
\begin{equation} \label{oooo}
\Pr  \Bigl( \max_{0 \leq s \leq t} X_s \geq r  \Bigr)
\leq (B t +x_0) r^{-1};
\end{equation}
\item[(ii)] for any $\eps>0$, a.s., for all but finitely many $t \in \Z^+$,
\[
\max_{0 \leq s \leq t} X_s  \leq   t (\log t)^{1+\eps}.
\]
\end{longlist}
\end{lemma}

We also state the following result on existence of passage-time moments
for one-dimensional stochastic processes, which is a simple
consequence of Theorem 1 of \cite{aim}.
%l2
\begin{lemma} \label{aimthm}
Let $(X_t)_{t \in \Z^+}$ be an $(\F_t)_{t\in\Z^+}$-adapted
stochastic process taking values in an unbounded
subset
$\mathcal{S}$ of $[0,\infty)$.
Suppose that $B>0$. Let
$\upsilon_B:= \min \{ t \in \N\dvt X_t \leq B \}$.
Suppose that
there exist $C \in (0,\infty)$ and  $\gamma \in [0,1)$ such that for any $t \in \Z^+$,
\[ \Exp [ X_{t+1}  - X_t    |  \F_t ] \leq
-C X_t ^\gamma   \qquad\mbox{on } \{ \upsilon_B > t \}.\]
Then, for any $p \in [0,1/(1-\gamma)]$ and  any $x \in \mathcal{S}$,
$\Exp [ \upsilon_B^p  |  X_0 = x]
< \infty$.
\end{lemma}

%s4.2 ###
\subsection{Exclusion-voter configurations}\label{config}

We introduce some more notation.
For $S \in \D \setminus\{ \D_0\}$ and $i \in \{1,\ldots,N\}$, let
%e1 ###
\begin{equation}\label{rt}
R_i:= R_i(S):= \sum_{j=1}^i n_j \quad\mbox{and}\quad T_i:= T_i(S):= \sum_{j=i}^N m_j.
\end{equation}

\begin{figure}

\includegraphics{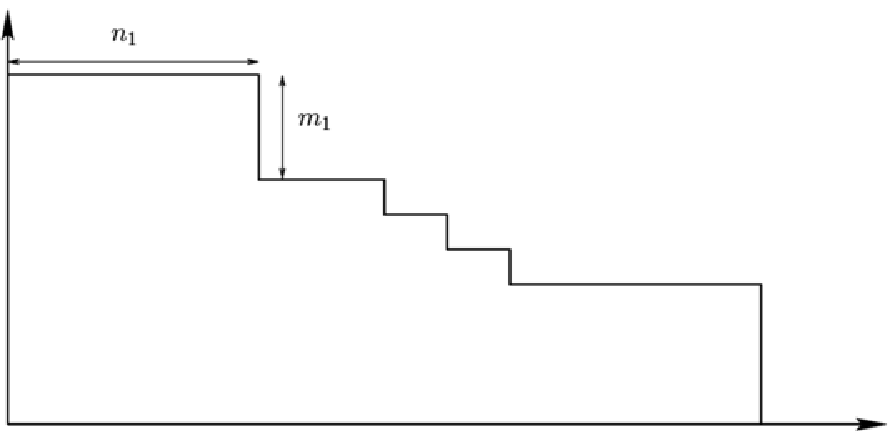}

\caption{An example of a staircase configuration.}\label{stair}
\end{figure}

It is convenient to represent
a configuration $S \in \D \setminus \{\D_0\}$
diagrammatically
as a right-down path in the quarter-lattice
$\Z^+ \times \Z^+$:
starting from
$(0, T_1)$, construct a walk by reading left-to-right
the configuration $S$ and, for each $0$ ($1$),
taking a unit step in the right (down)
direction. Thus, the walk  starts with
a step to the right and ends
at $(R_N,0),$ after $|S|$ steps.
See
Figure~\ref{stair} for the case of $S$
as given by (\ref{configex}).

%f2 ###

The lattice squares
of $\Z^+ \times \Z^+$
bounded by the
right-down path
determined by $S$
constitute a polygonal
region in the plane that
we call the \textit{staircase}
corresponding to $S$.
With this representation of the configuration space,
the exclusion-voter  model
can be viewed as a  growth/depletion
process on staircases.
For instance, exclusion
moves are particularly simple
in this context, corresponding
to adding  or removing  a square at a corner.

As well as $\D_0$,
we   introduce special notation for one more configuration. Set
%e2 ###
\begin{equation}
\label{d1def} \D_1: = \ldots 11101000\ldots,
\end{equation}
the configuration with $N(\D_1)=1$ and vector representation $(1,1)$.

We now introduce notation
for the changes in configuration
brought about by voter and exclusion
moves.
Given the staircase
of $S$, there are $2N+1$
`corners' representing $10$'s and $01$'s alternately, of which
$N+1$ are $10$'s and $N$ are $01$'s.
In the staircase representation, these corners
have coordinates $(R_i,T_{i+1})$,
$i \in \{0,\ldots,N\}$ (for $10$'s)
and $(R_i, T_i)$, $i \in \{1,\ldots,N\}$ (for $01$'s),
where $R_0=T_{N+1}=0$.
Enumerate the $10$'s
left-to-right in the configuration $S$
by $0,1,\ldots,N$, and similarly
the $01$'s by $1,\ldots,N$.

For $j\in \{0,1,\ldots,N\}$,
let $v_j^{10 \mapsto 00} (S)$ (resp.,
$v_j^{10 \mapsto 11} (S)$) denote the
configuration obtained from  $S$
by performing a voter move
changing the $j$th $10$  to $00$ (resp., $11$).
Similarly, for $j \in \{1,\ldots,N\}$,
let $v_j^{01 \mapsto 00} (S)$,
$v_j^{01 \mapsto 11} (S)$ denote the\vspace*{-1pt}
configuration obtained from the two
possible voter moves at  the $j$th $01$. We use\vspace*{-1pt}
analogous notation for exclusion moves:
$e_j^{10 \mapsto 01} (S)$ ($j \in \{0,\ldots,N\}$),
$e_j^{01 \mapsto 10} (S)$ ($j \in \{1,\ldots,N\}$).

To conclude this section, we sketch the (elementary)
proof of Proposition~\ref{prop1}.

\begin{pf*}{Proof of Proposition~\ref{prop1}}
It is not hard to see that
$\D_0$ is an absorbing state for the
pure voter model ($\beta =1$) and for the left-moving
totally asymmetric exclusion process ($\beta=0, p=1$),
hence also for the mixture model under
$\Pr_{\beta,1}$ for any $\beta \in [0,1]$.

To show that all states within $\D$ communicate,
it suffices to show that $\Pr_{\beta,p} ( \xi_{t+k} = S_1
 |  \xi_t = S_0 ) > 0$ for some $k\in \N$ for each of the following:
\begin{longlist}[(iii)]
\item[(i)] $S_0 = \D_0$, $S_1 = \D_1$;
\item[(ii)] $S_0 = \D_1$, $S_1 = \D_0$;
\item[(iii)]  any $S_0$ with $|S_0| \geq 2$ and some $S_1$ with $|S_1|=
|S_0| +1$;
\item[(iv)]   any $S_0$ with $|S_0| \geq 3$ and some $S_1$ with $|S_1| \leq
|S_0| -1$;
\item[(v)] any $S_0$ with $|S_0| \geq 3$ and any $S_1$,
where $S_1$ is identical to $S_0$ apart from
in a single position $j \in \{2,3,\ldots,|S_0|-1\}$.
\end{longlist}
In other words, given that moves of types (i)--(v) can occur,
it is possible (with positive probability)
to step, in a finite number of moves,
between any two configurations
in $\D$
by first adjusting the length of the configuration via moves
of types
(i)--(iv) and then flipping the states in the interior of the
configuration via moves of type (v).
Similarly, to show that all states in $\D \setminus \{ \D_0\}$
communicate, it suffices to show that all moves
of types (iii)--(v) have positive probability.

It is not hard to see that voter moves can perform
moves of types (ii), (iii) and (iv) in a single step (i.e., with $k=1$).
Similarly, exclusion moves with $p<1$ can perform
moves of types (i) and (iii) in one step, while
exclusion moves with $p>0$ can perform moves of
types (ii) and (iv), possibly needing multiple steps.
We claim that moves of type (v) can be performed
provided: (a)~$\beta \in (0,1)$;  or (b) $\beta =0$ and
$p \in (0,1)$.

In case (a),
suppose  that we need to replace a $0$ by a $1$ in the interior
of a given configuration. If $p<1$,  we may perform a voter
move on the first $10$ to the left of the position to be changed
and then, if necessary, perform successive $10 \mapsto 01$
exclusion moves to `step' the $1$ into the desired position. If $p>0$,
an analogous procedure works, starting from the first $01$ to the \textit{right}. On the other hand,
if we need to replace a $1$ by a $0$, a similar argument applies.

In case (b), we cannot use voter moves, but both types of
exclusion move are permitted, so we can `bring in' any $0$ ($1$)
from outside the disordered region, rearrange as necessary
and `take out' the excess $1$ ($0$) to the other boundary.

It follows that moves of types (ii)--(v) are possible,
provided $\beta \neq 1$ and $(\beta,p) \notin \{(0,0),(0,1)\}$,
and all (i)--(v) are possible if we additionally impose the condition $p<1$.

To complete the proof, we need to demonstrate aperiodicity
in the case where $\beta \neq 1$, $p<1$ and $(\beta,p) \neq (0,0)$,
where all states communicate. Since $\beta \neq 1$, exclusion
moves may occur.
Moreover, every configuration other than $\D_0$
contains at least one pair of each type ($01$ and $10$).
Hence, there is a positive probability
that a configuration other than $\D_0$ remains unchanged at a given step
(when a proposed exclusion move fails to occur).
Thus, since all states communicate,
we have aperiodicity.
\end{pf*}

%s4.3 ###
\subsection{Lyapunov function lemmas}\label{lemmas}

Throughout this paper, Lyapunov-type functions
will be primary tools. In this section,
we introduce some of our functions
and give some preliminary results.
Recall the definitions of $R_i$, $T_i$ from (\ref{rt}).
In \cite{bfmp},
much use was made of
the functions $f_1, f_2$  defined
for $S \in \D \setminus\{ \D_0 \}$ by
\[
f_1 (S):= \sum_{i=1}^N m_i R_i = \sum_{i=1}^N n_i T_i,\qquad
f_2(S):=  \frac{1}{2}  \Biggl( \sum_{i=1}^N m_i R_i^2 +
\sum_{i=1}^N n_i T_i^2  \Biggr),
\]
and by $f_1(\D_0)=f_2(\D_0)=0$.
Note that, with the diagrammatical representation
described in Section~\ref{config},
$f_1$ is the area of the staircase; for example,~for
$S$ given by (\ref{configex}), $f_1(S)=162$.

In the present paper, we introduce some more Lyapunov-type functions
that will prove valuable: these include $\rho^2$ (see (\ref{rhodef}) below),
$\phi_\alpha$ for $\alpha >0$ (see (\ref{phidef}) below)
and  $g$, which we define shortly.
First, we state some inequalities involving $f_1$ and $f_2$.
%l3
\begin{lemma}
For any $S \in \D$, we have
\begin{eqnarray}\label{ineq1}
&\frac{1}{2} |S| \leq f_1 (S) \leq \frac{1}{4} |S|^2 \quad\mbox{and}\quad
\frac{1}{4} |S|^2 \leq f_2 (S) \leq \frac{1}{8} |S|^3;&
\\\label{in1}
&f_2(S) \leq  |S| f_1 (S) \leq 2 (f_1(S))^2.&
\end{eqnarray}
\end{lemma}

\begin{pf}
The inequalities in (\ref{ineq1})  are  in
Lemma 4.1 of \cite{bfmp}. For (\ref{in1}), we have that for $S \in \D$,
%e3 ###
\begin{equation}\label{in3}
f_2 (S) \leq \frac{1}{2}  \Biggl( \sum_{i=1}^N m_i R_i + \sum_{i=1}^N n_i T_i  \Biggr)
\cdot ( R_N + T_1 ) = f_1 (S) \cdot |S|
\end{equation}
since, by (\ref{rt}),
$R_i \leq R_N$ and $T_i \leq T_1$ for $1 \leq i \leq N$.
Then, from (\ref{in3}) and the first $f_1$ inequality in~(\ref{ineq1}),
we obtain (\ref{in1}).
\end{pf}

The next lemma collects
formulae that we will need for the expected
increments of $f_1 (\xi_t)$ and $f_2 (\xi_t)$,
obtained from
(7.2), (5.3) and (6.3) in \cite{bfmp}.
Note that (\ref{0704a}) means that $f_2(\xi_t)$
is a martingale when $\beta=1$.
%l4
\begin{lemma}
If $S \in \D \setminus \{ \D_0\}$ and $\beta, p \in [0,1]$, then
\begin{eqnarray}\label{0704c}
&&\Exp_{\beta,p} [ f_1 (\xi_{t+1}) - f_1 (\xi_t)  |\xi_t = S ]\nonumber
\\[-8pt]\\[-8pt]
&&\quad= (1-\beta) \frac{N(1-2p)+(1-p)}{2N+1} - \beta \frac{N}{2N+1}; \nonumber\\\label{0704a}
&&\Exp_{\beta,p} [ f_2 (\xi_{t+1}) - f_2 (\xi_t)  |\xi_t = S ]\nonumber
\\[-8pt]\\[-8pt]
&&\quad = (1-\beta)  \Biggl( \frac{1}{2} + \frac{(1/2)-p}{2N+1} - \frac{2p-1}{2N+1} \sum_{i=1}^N (R_i + T_i)  \Biggr).\nonumber
\end{eqnarray}
\end{lemma}

Next, we define the  function $g$, which
captures most of $f_1$, in a sense made
precise in Lemma~\ref{lem0} below.
For $S \in \D \setminus \{\D_0\}$,
let $K = K(S)$ be the smallest
member of $\{1,\ldots,N\}$ for which
$R_K T_K = \max_{1 \leq k \leq N } \{ R_k T_k \}$.
Then, for $S \in \D \setminus \{\D_0\}$, set
\begin{equation}\label{xydef}
X(S):= R_K,\qquad   Y(S):= T_K
\end{equation}
and put $X(\D_0)=Y(\D_0)=0$. For $S \in \D$, we then define
\begin{equation}\label{gdef}
g (S):= X(S) Y(S) = \max_{1 \leq k \leq N} \{ R_k T_k \},
\end{equation}
where   $\max \varnothing:= 0$.
With the representation described in Section~\ref{config},
$g$ is the area of the largest rectangle
that can be inscribed in the staircase.
%l5
\begin{lemma}\label{lem0}
For any  $S \in \D \setminus\{ \D_0 \}$,
\begin{equation}\label{aa1}
f_1 (S) \geq g  (S) \geq \frac{ f_1(S)}{ 1+\log f_1 (S)}.
\end{equation}
\end{lemma}

\begin{pf}
We start with a geometrical argument that will
yield the stated results
via the staircase representation
of configurations $S$.
Define   $r_a(x):= a/x$ for $a>0$ and $x>0$. For $a>0$ and $b \geq 1$,
let $R(a,b)$ denote the region
defined by
\[
R(a,b):= \bigl\{ (x,y) \in \R^2\dvt  0 \leq x \leq b, 0 \leq y \leq (a/x) \1_{\{ x \geq 1\}} +a\1_{\{ x<1\}} \bigr\}.
\]
Then,
with $| \cdot |$ denoting  Lebesgue measure on $\R^2$,
\[
|R(a,b)| = a + \int_1^b (a/x) \,\ud x = a + a \log b.
\]
Let $h\dvtx [0,\infty) \to [0,c]$ be a non-increasing bounded function such that $h(x)=c$ for $0 \leq x <1$, $h(d)=0$
and $h(x) \geq 1$ for $0 \leq x < d$, where $c \geq 1$ and $d \geq 1$.
Define
\[
M:= M(h):= \{ (x,y) \in \R^2\dvt  0 \leq x \leq d, 0 \leq y \leq h(x) \}.
\]

Let $a_0:= \sup\{ a > 0\dvt  \{ r_a (x)\dvt x>0 \} \cap M \neq \varnothing \}$, that is,~the greatest
value of $a$ for which a curve $r_a(x)$ intersects the region $M$.
Then, let $x_0$ be such that $(x_0,r_{a_0}(x_0)) \in M$.
Let $B(M)$ denote
the rectangle
with vertices $(0,0), (x_0,0), (0, r_{a_0} (x_0))$ and
$(x_0, r_{a_0}(x_0))$;
then
$|B(M)|=x_0 (a_0/x_0) = a_0$. Moreover,
it is clear that $B(M) \subseteq M$ and
$M \subseteq R(a_0,d)$, so
\begin{equation}\label{pp0}
| B(M)| \leq |M| \leq | R(a_0,d)| = a_0 (1+ \log d).
\end{equation}
So, using the fact that $r_{a_0} (d) = a_0/d \geq 1$,
we obtain from (\ref{pp0}) that
\begin{equation}\label{pp}
1 \leq \frac{ |M|}{|B(M)|} \leq 1 +\log d \leq 1 + \log a_0 = 1 + \log |B(M)|.
\end{equation}

We now translate the above argument into a proof of the lemma.
Fix a configuration $S \in \D \setminus \{ \D_0\}$
with block representation $(n_1,m_1,\ldots,n_N,m_N)$. For $x \geq 0$, define
\[
j_S(x):= \max \Biggl\{ j \in \Z^+,  j \leq N\dvt  \sum_{i=1}^j n_i \leq x \Biggr\}\quad \mbox{and}\quad
h_S(x):= \sum_{i=j_S(x)+1}^N m_i,
\]
where we interpret an empty sum as zero. Set $c_S=\sum_{i=1}^N m_i$
and $d_S=\sum_{i=1}^N n_i$.
Then $h_S(x)=c_S$
when $0 \leq x <1$, since $n_1 \geq 1$. Also, $h_S(x)=0$ for $x \geq d_S$ and
$h_S(x) \geq m_N \geq 1$ for $0 \leq x < d_S$.
Therefore, $h_S$ is a function of the form of $h$ in the first paragraph of the present proof.
In particular, $|M(h_S)|=f_1(S)$ and
$|B(M(h_S))|=g(S)$. Thus, (\ref{pp}) implies
(\ref{aa1}).
\end{pf}

%s5 ###
\section{Non-existence of passage-time moments}\label{moments}

For $t \in \Z^+$, let $\F_t$ denote the $\sigma$-field generated by $(\xi_s; s \leq t)$.
Recall the definitions
of $X(S)$ and $Y(S)$ from (\ref{xydef}).
For convenience, we set $X_t:= X(\xi_t)$, $Y_t:= Y(\xi_t)$ and
consider the auxiliary $(\F_t)$-adapted
process $(\tilde \xi_t )_{ t \in \Z^+}$
defined by
$\tilde \xi_t:= (X_t,Y_t) = (X(\xi_t), Y(\xi_t) )$;
$\tilde \xi_t$ takes values in the quarter-lattice $\Z^+  \times \Z^+$
and
$\xi_t = \D_0$
if and only if $\tilde \xi_t = (0,0)$.
Let $\sigma_{x,y}$ be the time for $\xi_t$ to hit the ground state configuration $\D_0$
(equivalently, the time taken for $\tilde \xi_t$ to hit the origin  $(0,0)$)
given the $\F_0$-event $\{X(\xi_0)=x, Y(\xi_0)=y\}$.
The crucial ingredient in the proof of non-existence of
moments will be the following result.
%l6
\begin{lemma}\label{l2}
Suppose that $p \leq 1/2$, $\beta \in [0,1]$.
There then exist $\delta>0$, $\gamma>0$ such that for all $x,y \in \Z^+$,
%e4 ###
\begin{equation}\label{prob1}
\Pr_{\beta,p}   \bigl( \sigma_{x,y} \geq \delta (x^2 +y^2)  \bigr)\geq \gamma
\end{equation}
and
for all $S \in \D \setminus
\{ \D_0\}$,
%e5 ###
\begin{equation}\label{prob2}
\Pr_{\beta,p}   \biggl( \tau  \geq \delta \frac{f_1(S)}{1+\log f_1 (S)}   \Big|  \xi_0 = S   \biggr)\geq \gamma.
\end{equation}
\end{lemma}

Note that (\ref{prob2}) is close to Conjecture 7.2 in \cite{bfmp}.
The proof of Lemma~\ref{l2} will be carried out in stages.
The next result gives control
over the size of the disordered region in the mixture process
of  voter model with  \textit{symmetric
or recurrent} exclusion ($p \geq 1/2$).
%l7
\begin{lemma}\label{uu}
Suppose that $p \geq 1/2$ and $\beta \in [0,1]$. Then,  for all $t\in \N$,
\[
\Pr_{\beta,p}  \Bigl( \max_{0 \leq s \leq t} | \xi_s | \leq 2 \sqrt{10} t^{1/2}  \Bigr)
\geq 0.95 - \frac{f_2 (\xi_0)}{10 t}.
\]
\end{lemma}

\begin{pf}
For
$p \geq 1/2$ and $\beta \in [0,1]$,
we have, from (\ref{0704a}),
that
$f_2 (\xi_t)$ satisfies
\[
\Exp_{\beta,p} [ f_2 (\xi_{t+1}) - f_2 (\xi_t)  |  \xi_t = S ] \leq \tfrac{1}{2}
\]
for all $S \in \D$. Applying
Lemma~\ref{lem13}(i) to $f_2 (\xi_t)$ with
$r  =10 t$ and $B=1/2$, (\ref{oooo}) implies that
\[
\Pr_{\beta,p}  \Bigl( \max_{0 \leq s \leq t} f_2 (\xi_t) \leq 10 t  \Bigr) \geq
1- \frac{(t/2)+f_2 (\xi_0)}{10t} = 0.95 - \frac{f_2(\xi_0)}{10t}.
\]
Then, using the  fact that $|S| \leq 2 (f_2(S))^{1/2}$ for any $S \in\D$ (by (\ref{ineq1})),
we obtain the result.
\end{pf}

Suppose that
$\xi_0 = S_0 \in \D \setminus \{ \D_0\}$ with corresponding $\tilde \xi_0 =
(x_0,y_0) \in \Z^+   \times \Z^+$, that is, $X(S_0)=x_0$ and $Y(S_0)=y_0$.
In order to enable us to identify positions within a configuration
$S \in \D \setminus
\{ \D_0\}$,
enumerate the positions
in the hybrid zone left-to-right as $1, 2, \ldots, |S|$.

We now return to the voter
plus transient ($p \leq 1/2$) exclusion model and define an
auxiliary \textit{coloured} process as
follows.
Set
$H:=\sum_{i=1}^{K(S_0)} (n_i + m_i)$,
recalling the definition of $K(S_0)$ from just above (\ref{xydef});
then, position $H$ in $S_0$
is necessarily occupied
by a $0$ and position $H+1$ by a $1$.
We \textit{colour} the $x_0$ $0$'s that
occupy positions  in $\{1,2,\ldots,H\}$
and the $y_0$ $1$'s that
occupy positions in $\{ H+1,\ldots,|S_0|\}$.
All other particles  are uncoloured.
Intuitively,
coloured
particles can be thought of as `high energy'. Next, we will define the evolution of the  colouring
corresponding to the process $(\xi_t)_{t \in \Z^+}$. We emphasize that
the colouring is associated with the \textit{particles} (i.e., $1$'s and $0$'s)
rather than the \textit{sites}.

The colour dynamics is as follows.
Exclusion moves do not alter any colour so that particles retain their colour-state after an exclusion move.
Voter moves affect the colouring of particles only
if the modified pair consists of exactly one
coloured particle, in which case the colouring is changed
as follows.  In a pair $01$ or $10$, suppose that the $1$ is
coloured while the $0$ is not; a~voter move to pair $00$ produces two uncoloured particles,
while a move to pair $11$ produces two coloured particles. On the other hand, if, in an unlike pair,
the $0$ is coloured and the $1$ not, then a voter move to $00$ produces two coloured particles and to $11$ produces two uncoloured particles.
We note the following facts about the dynamics:
\begin{longlist}[(a)]
\item[(a)] uncoloured $1$'s remain to the left of any coloured $1$'s and
uncoloured $0$'s remain to the right of coloured $0$'s;
\item[(b)] a necessary condition   to be in the ground state configuration
$\D_0$ is that  the  set of
coloured particles consists only of a (possibly empty) block of coloured $1$'s at the left boundary
of the hybrid zone
and a (possibly empty) block of coloured $0$'s at the right boundary.
\end{longlist}

With $\xi_0=S_0 \in \D$, for $t \in\N$,
let $\xi^*_t$ denote the configuration
$\xi_t$ with the associated colouring
as determined by $(\xi_0,\ldots,\xi_t)$ according to the mechanism
just described.

Let $\F^*_t$ denote the $\sigma$-field generated by $(\xi^*_s; s \leq t)$.
Define the $\F^*_t$-measurable random
variables $\ell_t$ and $r_t$ as follows. Let $\ell_t$
be the position (measured from the left end of the hybrid zone)
of the leftmost coloured $1$ in $\xi^*_t$ and $r_t$ be the position of the rightmost coloured $0$ in $\xi_t^*$;
initially, $r_0+1=\ell_0$, by construction.

As the process evolves, coloured $1$'s may end up to the left of coloured
$0$'s. We define an auxiliary process $(\zeta_t)_{t \in \Z^+}$
to keep track of such configurations.
Informally, when $\ell_t < r_t$, $\zeta_t$ will be the portion of $\xi_t$ between positions
$\ell_t$ and $r_t$. More formally, we introduce a
holding state $\D_0^*$ and set $\zeta_t = \D_0^*$
if $\ell_t \geq r_t$. If $\ell_t < r_t$,
then the configuration $\xi^*_t$ induces
a finite string of $0$'s and $1$'s obtained
by extracting the segment of $\xi_t^*$
between positions $\ell_t$ and $r_t$ (inclusive);
we call this string $\zeta_t$.
Then $(\zeta_t)_{t \in \Z^+}$ is an $(\F_t^*)$-adapted process with $\zeta_0 = \D^*_0$.

Note that, when it is not in state $\D^*_0$,
$\zeta_t$ contains  only \textit{coloured} particles when colours are transposed from $\xi_t^*$. The idea now
is that when $p \leq 1/2$, $\zeta_t$ behaves
like the mixture of voter and $p \geq 1/2$ exclusion, except that the presence of uncoloured
particles in $\xi^*_t$ causes it to `slow down'; we therefore aim for a version
of Lemma~\ref{uu} in this case. This is the next result.
%l8
\begin{lemma}\label{uuu}
Suppose that $p \leq 1/2$ and $\beta \in [0,1]$.
Then, for all $t \in \N$,
\[
\Pr_{\beta,p}  \Bigl( \max_{0 \leq s \leq t} | \zeta_s | \leq 2 \sqrt{10}   t^{1/2}  \Bigr)\geq 0.95.
\]
\end{lemma}
\begin{pf} We compare the process $(\zeta_t)_{ t \in \Z^+}$ to an independent
copy $\xi' = (\xi'_t)_{ t \in \Z^+}$ of the process~$\xi$.
We define $f^*_2(\zeta_t)$ analogously
to $f_2(\xi_t)$, but counting only the (coloured)
particles in region $\zeta_t$, that is, coloured $1$'s to the
left of coloured $0$'s and coloured $0$'s to the right of
coloured $1$'s.
Suppose
that we were to permit the initial configuration
$\xi^*_0 = \D'_0:=\ldots 000111 \ldots,$
where all $0$'s and $1$'s are coloured,
so that $\zeta_0=\D_0^*$. Then, by
a simple reflection argument,
the process $(\zeta_t)_{t \in \Z^+}$ embedded
in $(\xi_t^*)_{t \in \Z^+}$ started from $\xi^*_0=\D'_0$
has the same distribution
under $\Pr_{\beta,p}$ as the process $(\xi_t)_{ t \in \Z^+}$ under $\Pr_{\beta,1-p}$
with initial state $\D_0$. So, in particular, Lemma~\ref{uu} holds with $\zeta_t$ instead of $\xi_t$,
given the initial configuration $\D'_0$;  using
the fact that $f^*_2(\zeta_0)=0$,
we then obtain the claimed result in this case.

Now, the presence of
uncoloured $1$'s to the left of coloured $1$'s
or uncoloured $0$'s to the right of coloured $0$'s
restricts the growth of $|\zeta_t|$;
hence, the claimed result also holds for
any permissible initial configuration for $\xi_0^*$ other than $\D'_0$.
(One can argue  rigorously by stochastic domination at this point.)
\end{pf}

\begin{pf*}{Proof of Lemma~\ref{l2}}
We first prove the statement (\ref{prob1}).
Let $\chi_t:=\chi(\xi^*_t)$ denote the number of coloured particles in $\xi^*_t$. Then
$( \chi_t )_{ t \in \Z^+}$ is $(\F^*_t)$-adapted and $\chi_0=\chi(\xi^*_0)=x_0+y_0$.
Also,
given $\chi_t = n$ for $n \in \N$, we have that $\chi_{t+1}=n$,
unless a voter move is performed on a pair with exactly one particle coloured,
in which case $\chi_{t+1}$ takes values $n-1,n+1$ with equal probability. Also,
if $\chi_t=0$, then $\chi_{t+1}=0$ as well. Thus, $\chi_t$ is a non-negative $(\F^*_t)$-martingale
with uniformly bounded jumps. It follows from Doob's submartingale
inequality applied to the non-negative submartingale $(\chi_t-(x_0+y_0))^2$, using the
fact that $\Exp [ (\chi_t-(x_0+y_0))^2 ] \leq t$, by the orthogonality of martingale
increments, that for any $z>0$,
\[
\Pr_{\beta,p}   \Bigl( \max_{0 \leq s \leq t} | \chi_s - (x_0+y_0)  |^2 \geq z  \Bigr)
\leq t/z.
\]
In particular, taking $z=100t$, this implies that for any $t \in \Z^+$,
\[
\Pr_{\beta,p}   \Bigl( \min_{0 \leq s \leq t} \chi_s \geq (x_0+y_0) - 10 t^{1/2}  \Bigr)
\geq 0.99.
\]
Taking $t=\delta^2 (x_0^2 + y^2_0)$ for some $\delta>0$, combining
the last display with Lemma~\ref{uuu}, we have that, with
probability at least $0.94,$ the two events
\begin{eqnarray*}
&&\Bigl\{ \min_{0
\leq s \leq t} \chi_s \geq (x_0+y_0) - 10 \delta (x_0^2 +
y^2_0)^{1/2} \geq (1-10\delta) (x_0 +y_0)  \Bigr\}
\quad\mbox{and }\\
&&\Bigl\{ \max_{0 \leq s \leq t} | \zeta_s | \leq 2 \sqrt{10} \delta (x_0^2+ y_0^2)^{1/2} \leq
2 \sqrt{10} \delta (x_0 + y_0)  \Bigr\}
\end{eqnarray*}
both occur (noting that $(x_0^2 + y_0^2)^{1/2} \leq (x_0 + y_0)$).
Choose $\delta$ small, say $\delta =0.01$. Then, with probability  at least $0.94$, the
total number $\chi_s$
of coloured particles up to time $t$ remains greater than $0.9 (x_0+y_0)$,
while the central overlap region $\zeta_s$
of coloured particles between the leftmost coloured $1$
and the rightmost coloured $0$
remains shorter than $0.1 (x_0+y_0)$.
Hence, there must remain at least one coloured $0$ to the left of any coloured $1$
or one coloured $1$ to the right of any coloured $0$.
By observation (b) above,
this excludes the possibility of $\xi_s = \D_0$ for any $s \leq t$,
where $t = \delta^2 (x^2_0 +y^2_0)$.
We thus obtain (\ref{prob1}).

To derive (\ref{prob2}), we use the fact that
$x_0^2+y_0^2 \geq 2 x_0 y_0 = 2 g (\xi_0)$
and then use (\ref{aa1}).
\end{pf*}

We are now nearly ready to complete the proofs of Theorems~\ref{thm1} and~\ref{thm2}.
The proofs proceed in a similar
way to the proof of Theorem 6.1 in \cite{bfmp}.

\begin{pf*}{Proof of Theorem~\ref{thm1}}
Suppose that $p \leq 1/2$.
Take $S_0 \in \D \setminus\{ \D_0 \}$.
Suppose, for the purpose of deriving a contradiction, that
$\Exp_{\beta,p} [ \tau^{1+\eps}  |  \xi_0 = S_0 ] < \infty$ for some $\eps>0$.
Let $\xi'=(\xi'_t)_{t \in \Z^+}$ be an independent copy of $\xi$ and $\tau'$
be the corresponding independent copy of $\tau$.
For any $t \in \Z^+$,  using the Markov property,
we obtain
\begin{equation}\label{0420a}
\Exp_{\beta,p} [ \tau^{1+\eps}  |  \xi_0 = S_0 ]
\geq
\Exp_{\beta,p}  \bigl[ \Exp_{\beta,p}  [
( t+\tau' )^{1+\eps}  |  \xi'_0 = \xi_t  ] \1_{\{\tau \geq t\}}   |  \xi_0 = S_0  \bigr].
\end{equation}
For the inner expectation in the  expression on the right-hand side
of (\ref{0420a}), we have, by
(\ref{prob2}), that there exist $\delta>0$, $\gamma>0$
such that for any $t \in \Z^+$,
\[
\Exp_{\beta,p}  [ (t+ \tau')^{1+\eps}  |  \xi_0' = \xi_t  ]
\geq \gamma  \biggl( \delta \frac{f_1 (\xi_t)}{1 + \log f_1 (\xi_t)}  \biggr)^{1+\eps}.
\]
Since, for any $\eps>0$,  $x^\eps > 1 + \log x$ for all $x$
sufficiently large,
there exists $\gamma' >0$ for which
%e6 ###
\begin{equation}\label{0420b}
\Exp_{\beta,p}  [ (t+ \tau')^{1+\eps}  |  \xi_0' = \xi_t  ]
\geq \gamma'  \bigl( f_1 (\xi_t)^{1-(\eps/2)}  \bigr)^{1+\eps}
\end{equation}
for any $t \in \Z^+$.
It follows from (\ref{0420a}) with (\ref{0420b})
that for some $\eps' \in (0,\eps)$ and some $C \in (0,\infty)$,
\[
\Exp_{\beta,p} [ \tau^{1+\eps} |  \xi_0 = S_0 ]
\geq  C \Exp_{\beta,p} \bigl[ (f_1 (\xi_t))^{1+\eps'} \1_{\{\tau \geq t\}}  |  \xi_0 = S_0 \bigr] = C \Exp_{\beta,p}
[  (f_1 (\xi_{t \wedge \tau} ))^{1+\eps'}  |  \xi_0 = S_0  ]
\]
for any $t \in \Z^+$,
using the fact that $f_1(\xi_{\tau} ) = f_1(\D_0) = 0$ a.s.
That is, given $\xi_0 = S_0$,
$(f_1 (\xi_{t \wedge \tau}))^{1+\eps'}$ is uniformly
bounded in $L^1$.

Hence, the assumption that $\Exp_{\beta,p} [ \tau^{1+\eps}  |  \xi_0 =S_0] < \infty$
implies that on $\xi_0 = S_0$,
the process $f_1 (\xi_{t \wedge \tau})$ is uniformly integrable and, trivially,
that $\tau < \infty$ a.s.; thus,
as $t \to \infty$,
$\Exp_{\beta,p} [ f_1 (\xi_{t \wedge \tau})  |  \xi_0 = S_0 ] \to \Exp_{\beta,p} [ f_1 (\xi_{\tau})   |  \xi_0 = S_0 ] =
f_1 (\D_0) = 0$.
However,
for $p \leq  1/2$ and $\beta \leq (1-2p)/(2-2p)$, it follows
from (\ref{0704c})  that   for any
$t \in \Z^+$ and any $S \in \D$,
\[ \Exp_{\beta,p} [ f_1 (\xi_{t+1} ) - f_1 (\xi_t)  |  \xi_t = S ] \geq 0.\]
By the submartingale property,
we then have that for all $t \in \Z^+$,
$\Exp_{\beta,p} [ f_1 (\xi_{t \wedge \tau})   |  \xi_0 = S_0 ] \geq   f_1 (S_0)
> 0$. We thus have the desired contradiction.
\end{pf*}

\begin{pf*}{Proof of Theorem~\ref{thm2}}
Suppose
that $S_0 \in \D \setminus \{ \D_0\}$ and, for a contradiction,
that
$\Exp_{\beta,p} [  \tau ^{2+\eps}  |  \xi_0 = S_0] < \infty$ for some $\eps>0$.
Then, for any $t \in \Z^+$,  similarly
to the proof of Theorem~\ref{thm1},
\[
\Exp_{\beta,p} [  \tau ^{2+\eps}  |  \xi_0 =S_0 ]  \geq
\Exp_{\beta,p}  \bigl[ \Exp_{\beta,p}  [
( t+\tau' )^{2+\eps}  |  \xi'_0 = \xi_t  ] \1_{\{\tau  \geq t\}}
 |  \xi_0 = S_0  \bigr].
 \]
Hence, for $p \leq 1/2$, using (\ref{prob2}), there exist $\gamma, \delta, \eps', \eps'' >0$ such that
\begin{eqnarray*}
\Exp_{\beta,p} [ \tau^{2+\eps}  |  \xi_0 = S_0 ]
& \geq&    \gamma \Exp_{\beta,p}  \bigl[ \bigl(t+\delta (f_1 (\xi_t ))^{1-(\eps/3)} \bigr)^{2+\eps}
\1_{\{\tau  \geq t\}}  |  \xi_0 = S_0  \bigr]\\
&  \geq &   C \Exp_{\beta,p} [ ( f_1 (\xi_{t \wedge \tau }))^{2+\eps'}  |  \xi_0 = S_0]
\geq
C   \Exp_{\beta,p} [ (f_2 (\xi_{t \wedge \tau }))^{1+\eps''}   |  \xi_0 = S_0],
\end{eqnarray*}
using (\ref{in1}) for the last inequality.
Hence, the process $f_2 (\xi_{t \wedge \tau })$ is uniformly integrable and thus,
as $t \to \infty$,
$\Exp_{\beta,p} [ f_2 (\xi_{t \wedge \tau })  |  \xi_0 = S_0 ]
\to \Exp_{\beta,p} [ f_2 (\xi_{\tau })   |  \xi_0 = S_0 ] =
f_2(\D_0)=0$.

However,
for $p \leq 1/2$ and $\beta \in [0,1]$, for
all $S \in \D$ and all
$t \in \Z^+$, it follows from (\ref{0704a})  that
\[ \Exp_{\beta,p} [ f_2 (\xi_{t+1} ) - f_2 (\xi_t)  |  \xi_t = S ] \geq 0.\]
Hence,
for all $t \in \Z^+$,
$\Exp_{\beta,p} [ f_2 (\xi_{t \wedge \tau })  |  \xi_0 = S_0 ]
\geq   f_2 (S_0)
> 0$, giving a contradiction.
\end{pf*}

%s6 ###
\section{Recurrence}\label{recprf}

We consider a new Lyapunov-type function
that generalizes $f_1$. For $\alpha \geq 0$,
set $\phi_\alpha (\D_0):= 0$ and
for $S \in \D \setminus \{ \D_0\}$, set
%e7 ###
\begin{equation}\label{phidef}
\phi_\alpha (S):= \sum_{i=1}^N \sum_{j=R_{i-1}+1}^{R_i}
\sum_{k=1}^{T_i} \frac{1}{(j+k)^\alpha} =
\sum_{i=1}^N \sum_{j=T_{i+1}+1}^{T_i}
\sum_{k=1}^{R_i} \frac{1}{(j+k)^\alpha};
\end{equation}
here, and
throughout this section,
we use the conventions $R_0:= 0$, $R_{N+1}:= R_N$,
$T_0:= T_1, T_{N+1}:= 0$. In particular, it
follows from (\ref{phidef}) that when $\alpha =0$,
$\phi_0 (S) = \sum_{i=1}^N n_i T_i = \sum_{i=1}^N m_i R_i = f_1 (S)$.
For convenience, we introduce the notation
\[
a_i(j):= (T_j +R_j + i)^{-\alpha}  \quad\mbox{and}\quad
b_i(j):= (T_{j+1} +R_j + i)^{-\alpha}.
\]
The next lemma gives an expression for the expected increments of $\phi_\alpha$.
%l9
\begin{lemma}\label{change}
Let $\beta \in [0,1]$ and $p\in [0,1]$.
Then, for any $S \in \D \setminus \{\D_0\}$ and any $t \in \Z^+$,
\begin{eqnarray}\label{phi1}
&&\Exp_{\beta,p} [ \phi_\alpha (\xi_{t+1}) - \phi_\alpha (\xi_t)
 |  \xi_t = S ]\nonumber
\\
&&\quad=   \frac{1-\beta}{2N+1}
\Biggl\{ -p \sum_{j=1}^N a_0(j)
+(1-p) \sum_{j=0}^N b_2(j)
\Biggr\}
\\
&&\qquad{}+ \frac{\beta}{2N+1}  \Biggl\{ N \sum_{j=1}^N
\bigl( a_1(j) - b_1 (j-1)
\bigr)
- \frac{1}{2} \sum_{j=2}^N a_1 (j)
- \frac{N+1}{2} b_1(N)  \Biggr\}
.\nonumber
\end{eqnarray}
\end{lemma}

\begin{pf}
Recalling the notation of Section~\ref{config},
we  write
\begin{eqnarray*}
D_j^{\mathrm{v,10}} (S)
&:=&   \phi_\alpha ( v_j^{10 \mapsto 00} (S) )
+   \phi_\alpha
( v_j^{10 \mapsto 11} (S) )  - 2\phi_\alpha (S) \qquad (j \in \{0,\ldots,N\}), \\
D_j^{\mathrm{v,01}} (S)
&:=&  \phi_\alpha ( v_j^{01 \mapsto 00} (S) )
+   \phi_\alpha
( v_j^{01 \mapsto 11} (S) )  - 2\phi_\alpha (S) \qquad (j \in \{1,\ldots,N\}), \\
D_j^{\mathrm{e,10}} (S)
&:=&   \phi_\alpha ( e_j^{10 \mapsto 01} (S) )
- \phi_\alpha (S) \qquad (j \in \{0,\ldots,N\}), \\
D_j^{\mathrm{e,01}} (S)
&:=&   \phi_\alpha ( e_j^{01 \mapsto 10} (S) )
- \phi_\alpha (S) \qquad (j \in \{1,\ldots,N\}).
\end{eqnarray*}
Summing over all
possible moves,
we have that
\begin{eqnarray}\label{aaq1}
&&\Exp_{\beta,p} [ \phi_\alpha (\xi_{t+1}) - \phi_\alpha (\xi_t)
 |  \xi_t = S ]\nonumber
 \\
 &&\quad=
\frac{\beta}{2N+1}  \Biggl\{
\frac{1}{2} \sum_{j=1}^N D_j^{\mathrm{v,01}} (S) + \frac{1}{2} \sum_{j=0}^N D_j^{\mathrm{v,10}} (S)
\Biggr\} \\
&&\qquad{}+ \frac{1-\beta}{2N+1}  \Biggl\{ -p \sum_{j=1}^N
D_j^{\mathrm{e,01}} (S) + (1-p) \sum_{j=0}^N
D_j^{\mathrm{e,10}} (S)  \Biggr\}.\nonumber
\end{eqnarray}

We now calculate
expressions for the terms in (\ref{aaq1}).
The reader might find it helpful here to
refer to a picture such as Figure~\ref{stair} in Section~\ref{config}.
We have that for $j \in \{1,\ldots,N\}$,
\[
D_j^{\mathrm{v,01}} (S) =
\sum_{i=1}^N  \bigl( a_1(i)
- b_1(i-1)  \bigr)
- a_1(j)
- a_1(j+1).\]
Also, for $j \in \{0,1,\ldots,N\}$, we have that $D_j^{\mathrm{v,10}} (S)$ is given by
\[
\sum_{i=1}^N   \bigl( a_1(i)
- \1_{\{ i \leq j\}} b_1(i-1)
- \1_{\{ i > j\}} b_1(i)   \bigr)
=
\sum_{i=1}^N   \bigl( a_1(i) - b_1(i-1)  \bigr)
+ b_1(j)
- b_1(N).
\]
Taking the computations for  $D_j^{\mathrm{v,01}} (S)$,  $D_j^{\mathrm{v,10}} (S)$ and summing, we have
\begin{eqnarray*}
&& \frac{1}{2} \sum_{j=1}^N D_j^{\mathrm{v,01}} (S) + \frac{1}{2} \sum_{j=0}^N D_j^{\mathrm{v,10}} (S)
\\
&&\quad=   \frac{2N+1}{2} \sum_{j=1}^N
\bigl( a_1(j) - b_1(j-1)  \bigr)
- \frac{1}{2} \sum_{j=1}^N \bigl( a_1(j)  + a_1(j+1) \bigr)
\\
&&\qquad{}+
\frac{1}{2} \sum_{j=0}^N b_1(j)
- \frac{N+1}{2} b_1(N) \\
&&\quad=  \frac{2N+1}{2} \sum_{j=1}^N
\bigl( a_1(j)
- b_1(j-1)  \bigr)
- \sum_{j=1}^N a_1(j)
+ \frac{1}{2} a_1(1)
\\
&&\qquad{}+ \frac{1}{2} \sum_{j=1}^N b_1(j-1)
- \frac{N+1}{2} b_1(N) \\
&&\quad=
N \sum_{j=1}^N
\bigl( a_1(j)
- b_1(j-1) \bigr)
- \frac{1}{2} \sum_{j=2}^N a_1(j)
- \frac{N+1}{2} b_1(N).
\end{eqnarray*}
For the (simpler) exclusion moves, we obtain
$D_j^{\mathrm{e,01}} (S) = -a_0(j)$ and
$D_j^{\mathrm{e,10}} (S) = b_2(j)$. Then,
combining all the computations, from  (\ref{aaq1}),
we obtain (\ref{phi1}).
\end{pf}

For the\vspace*{2pt} rest of this section, we will be interested in the properties
of $\phi_1$.
%l10
\begin{lemma}\label{phibnd}
For any $S \in \D \setminus \{ \D_0\}$, we have $\phi_1 (S) \geq   \log (|S|/4)$.
\end{lemma}
\begin{pf} Suppose that $S \in \D \setminus \{ \D_0\}$. From (\ref{phidef}),
we have that
\[ \phi_1 (S) \geq \sum_{i=1}^N \sum_{j=R_{i-1}+1}^{R_i} \frac{1}{1+j}
= \sum_{j=1}^{R_N} \frac{1}{j+1}   \geq \int_1^{R_N} \frac{ \ud x}{1+x} \geq \log ( R_N /2),\]
using monotonicity for the second inequality.
Similarly, (\ref{phidef}) gives
$\phi_1(S) \geq \log (T_1/2)$. Thus,
$\phi_1(S) \geq \log ( \max \{ R_N, T_1 \} /2 )$, which yields the result.
\end{pf}

The following lemma is the key to this section.
%l11
\begin{lemma}\label{nnn3}
Suppose that $\beta \geq 4/7$ and $p \in [0,1]$.  Then,
for any $S \in \D \setminus\{ \D_0\}$,
\[
\Exp_{\beta,p} [ \phi_1 (\xi_{t+1}) - \phi_1 (\xi_t)
 |  \xi_t = S ] \leq 0.\]
\end{lemma}

\begin{pf}
For ease of notation during this proof, set
$\Delta(S):= \Exp_{\beta,p} [ \phi_1 (\xi_{t+1}) - \phi_1 (\xi_t)
 |  \xi_t = S ]$.
It is clear from (\ref{phi1}) that
$\Delta(S)$ is non-increasing in $p$
and so it suffices to consider the
case $p=0$.
(\ref{phi1}) then implies that in this case, $\Delta(S)$ is given by
\[
\frac{\beta}{2N+1}  \Biggl\{ N \sum_{j=1}^N
\bigl( a_1(j) - b_1(j-1)  \bigr)
- \frac{1}{2} \sum_{j=2}^N a_1(j)
- \frac{N+1}{2} b_1(N)  \Biggr\}
+ \frac{1-\beta}{2N+1}
\sum_{j=0}^N b_2(j).
\]
We rewrite this last expression
by setting
$\gamma:= (1-\beta)/\beta \in [0,\infty)$ to
obtain
\begin{eqnarray}\label{0407a}
\frac{2N+1}{\beta} \Delta (S)
&=&   N \sum_{j=1}^N  \bigl(a_1(j) - b_1 (j-1)  \bigr)\nonumber
\\[-8pt]\\[-8pt]
&&{}-\frac12 \sum_{j=2}^N a_1(j) -\frac 12 \frac{N+1}{R_N+1}
+\gamma \sum_{j=1}^{N+1} b_2(j-1).\nonumber
\end{eqnarray}
We need to show that the right-hand side of (\ref{0407a})
is non-positive. Since this quantity is non-decreasing
in $\gamma$, it suffices
to consider the case $\gamma =3/4$, corresponding
to $\beta = 4/7$.
Set
\begin{eqnarray*}
\tilde \Delta (S):=
N\sum_{j=1}^N \bigl( a_1(j) - b_1(j-1)  \bigr)- \frac{1}{2} \sum_{j=2}^N
a_1(j) -\frac 12 \frac{N+1}{R_N+1}
+ \frac{3}{4}
\sum_{j=1}^{N+1} b_1(j-1)
\end{eqnarray*}
so that, from (\ref{0407a}),
$\Delta (S) \leq \frac{\beta}{2N+1}\tilde \Delta(S)$ since $b_1(j) \geq b_2(j)$.

Write $A_N:=1+m_1+m_2+\cdots+m_N$,
$D_0:=0$ and, for $i \in \{1,\ldots,N\}$,
$D_i:=(n_1-m_1)+\cdots+(n_i-m_i)$
so that $R_{j-1}+T_j+1=A_N+D_{j-1}$.
We then have that
\begin{eqnarray}\label{eq1}
\tilde\Delta (S)
&= &
N\sum_{j=1}^N \biggl(\frac1{A_N+D_{j-1}+n_j}-\frac1{A_N+D_{j-1}} \biggr)-
\sum_{j=2}^N \frac{1/2}{A_N+D_{j-1}+n_j} \nonumber\\[-8pt]\\[-8pt]
&&{}  - \frac{(N+1)/2}{A_N+D_N}+
\frac{3}{4} \sum_{j=1}^{N+1}\frac1{A_N+D_{j-1}}
=
\frac {1/2}{A_N+n_1}+ H_N(S),\nonumber
\end{eqnarray}
where we have introduced the notation,
for $k \in \{1,\ldots,N\}$,
\[ H_k (S):=  \sum_{j=1}^k
\biggl(\frac{N-1/2}{A_N+D_{j-1}+n_j} -
\frac{N-3/4}{A_N+D_{j-1}} \biggr)
-\frac{(N+k-1)/4}{A_N+D_k}. \]
We now claim that if $N \geq 2$, then for any $k \in \{2,\ldots,N\}$,
%e8 ###
\begin{equation}\label{hclaim}
H_k (S) \leq H_{k-1} (S).
\end{equation}
We then have, from (\ref{eq1}) and (\ref{hclaim}), that
for $N \geq 2$,
\begin{eqnarray*}
\tilde \Delta (S) & =&    \frac {1/2}{A_N+n_1}+ H_N (S)
\leq \frac {1/2}{A_N+n_1}+ H_1 (S)
= \frac{N}{A_N +n_1} -
\frac{N-3/4}{A_N}
-\frac{N/4}{A_N+n_1 - m_1}  \\
& \leq&
\frac{N}{A_N +n_1} -
\frac{N-3/4}{A_N+n_1}
-\frac{N/4}{A_N+n_1 } = \frac{(3-N)/4}{A_N+n_1}.
\end{eqnarray*}
Thus, $\tilde\Delta (S) \leq 0$ and hence
$\Delta(S) \leq 0$ also,
for all $S \in \D \setminus \{ \D_0\}$
with $N(S)\ge 3$.

Let us now verify the claim (\ref{hclaim}).
We have that for $k \in \{ 2, \ldots, N \}$,
\begin{eqnarray*}
H_k (S) &=& \sum_{j=1}^{k-1}
\biggl(\frac{N-1/2}{A_N+D_{j-1}+n_j} -
\frac{N-3/4}{A_N+D_{j-1}} \biggr)
+ \frac{N-(1/2)}{A_N + D_{k-1} + n_k} \\
&&{}- \frac{N-(3/4)}{A_N + D_{k-1}}
- \frac{(N+k-1)/4}{A_N + D_k} \\
&=& \sum_{j=1}^{k-1}
\biggl(\frac{N-1/2}{A_N+D_{j-1}+n_j} -
\frac{N-3/4}{A_N+D_{j-1}} \biggr)
+  \biggl[ \frac{(3N-k-1)/4}{A_N + D_{k-1} + n_k}
- \frac{N-(3/4)}{A_N + D_{k-1}} \biggr] \\
&&{}+\biggl[ \frac{(N+k-1)/4}{A_N + D_{k-1} + n_k} - \frac{(N+k-1)/4}{A_N + D_k}  \biggr],
\end{eqnarray*}
where we have split the term with the denominator $A_N +D_{k-1}+n_k$
into two parts. Note that for all $j$ we have $A_N+D_{j-1}+n_j \geq A_N+D_{j-1}$ and
also $A_N+D_{j-1}+n_j=A_N+D_{j}+m_j \geq A_N+D_{j}$. Therefore,
applying these inequalities separately
to the two terms in square brackets in the last display,
we verify the claim (\ref{hclaim}) since
\begin{eqnarray*}
H_k(S) \leq
\sum_{j=1}^{k-1}
\biggl(\frac{N-1/2}{A_N+D_{j-1}+n_j} -
\frac{N-3/4}{A_N+D_{j-1}} \biggr)
+  \biggl[ \frac{(-N-k+2)/4}{A_N + D_{k-1}}  \biggr]
= H_{k-1}(S).
\end{eqnarray*}

To complete the proof,
we show that $\Delta(S) \leq 0$ for $N(S) \in \{1,2\}$ also.
For $N=1$, writing the right-hand side of the $\beta=4/7$ case of (\ref{0407a}) over a common
denominator, we have
\begin{eqnarray*}
\frac{21}{4}\Delta (S)  =
-\frac{n_1 m_1 (n_1-m_1)^2+ 13(n_1+m_1)+4(1+n_1m_1)+14(n_1^2+m_1^2)+5(n_1^3+m_1^3)}{4(1+ m_1+n_1)(1+ m_1)(1+ n_1)(2+ m_1)(2+ n_1)}
,
\end{eqnarray*}
which is negative.
Finally, for $N=2$,  from (\ref{0407a}) and some tedious algebra,
we obtain $35 \Delta (S) / 4= - Q/R$,
where $R=4(m_1+m_2+n_1+1)(m_1+m_2+1)(m_2+n_1+n_2+1)(m_2+n_1+1)(n_1+n_2+1)(m_1+m_2+2)(m_2+n_1+2)(n_1+n_2+2)$
and
\[
Q=m_1 n_1^4 (2m_1^2+5n_1^2-m_1 n_1)+ n_2 m_2^4
(2n_2^2+5m_2^2-n_2 m_2)+\mbox{$244$ positive terms,}
\]
as   can be readily checked in \texttt{Maple}, for instance. Since $2x^2+5y^2-xy$
is always non-negative, we conclude that $\Delta(S)\le 0$ in this last case also.
\end{pf}

\begin{pf*}{Proof of Theorem~\ref{recthm}}
Lemma~\ref{nnn3}
shows that for $\beta \geq 4/7$,
$(\phi_1(\xi_t))_{t \in \Z^+}$
is a
supermartingale on $\xi_t \in \D \setminus \{ \D_0 \}$.
Since, by Lemma~\ref{phibnd},
$\phi_1(S) \to \infty$ as $| S| \to \infty$,
we can use Theorem 2.2.1 of \cite{fmm}
to complete the proof of the theorem.
\end{pf*}

%s7 ###
\section{Existence of passage-time moments}\label{secex}

Our main tool in this section will be Lemma~\ref{aimthm} applied with
the Lyapunov function $f_2$.  Our first result is a bound on the
expected increments of $f_2$.
%l12
\begin{lemma}\label{f2lem}
Suppose that  $\beta \in [0,1)$ and $p>1/2$. There then
exists $C \in (0,\infty)$ such that for
all but finitely many $S \in \D$,
\[ \Exp_{\beta,p} [ f_2 (\xi_{t+1}) - f_2 (\xi_{t})
 |  \xi_t = S ] \leq - C ( f_2 (S) )^{1/6}.\]
\end{lemma}
\begin{pf}
It follows directly from (\ref{rt}) that
$R_N + T_1 = |S|$ and $R_i+T_i \geq N$ so that
\begin{equation}\label{0825a}
\sum_{i=1}^N (R_i + T_i ) \geq \max \{ |S|, N^2 \} \geq N |S|^{1/2}.
\end{equation}
We see from (\ref{0704a}) with (\ref{0825a}) that
for $p>1/2$, $\Exp_{\beta,p} [ f_2 (\xi_{t+1}) - f_2 (\xi_{t})
 |  \xi_t = S ]$ is at most
\[  \frac{1-\beta}{2} - (1-\beta) (2p-1) \frac{N
|S|^{1/2} }{2N+1} \leq \frac{1-\beta}{2} - \frac{1}{3}(1-\beta)
(2p-1)  |S|^{1/2} \]
since $N \geq 1$. The result follows since $|S|^{1/2} \geq f_2(S)^{1/6}$,
from (\ref{ineq1}).
\end{pf}

\begin{pf*}{Proof of Theorem~\ref{thmex2}}
The $\beta = 1$ case of the theorem follows from
Theorem~\ref{bfmpthm}(ii). Now, suppose that $\beta \in [0,1)$.
Applying Lemma~\ref{aimthm} with $X_t = f_2 (\xi_t)$
and using Lemma~\ref{f2lem} shows that the
hitting time of a finite subset of $\D$ has finite $(6/5)$th moment.
Since $\D_0$ is accessible from any state, it follows that
$\tau$ also has finite $(6/5)$th moment.
\end{pf*}
%r
\begin{remark*}
The exponent $1/6$ in Lemma~\ref{f2lem} may not be the best possible. However,
(\ref{0704a}) applied to
the configuration $n_1=n_2=\cdots=n_{N-1}=1,n_N=N^2$ and
$m_1=N^2, m_2 = \cdots = m_N = 1$ shows that one cannot increase the
exponent to more than $1/4$ in general. Hence, the method used in this
section seems unable to
prove existence of moments greater than $4/3$; see the remark immediately following the statement
of Theorem~\ref{thmex2}.
\end{remark*}

To prove Theorem~\ref{thmex1}, we will again apply Lemma~\ref{aimthm},
but this time we will take $X_t = f_2 (\xi_t)^M$ for arbitrary $M \in [1,\infty)$.
To study the increments of this process, we recall some facts about $f_2$
under exclusion moves; compare (5.1) and (5.2) in \cite{bfmp}.
We have that
\begin{eqnarray}\label{0705a}
f_2 ( e_j^{10 \mapsto 01} (S) ) &=& f_2 (S) + 1 + R_j + T_{j+1}\qquad  (j
\in \{0,\ldots,N\}),\\
\label{0705b}
f_2 ( e_j^{01 \mapsto 10} (S) ) &=& f_2 (S) + 1 - R_j - T_{j}\qquad  (j
\in \{1,\ldots,N\}),
\end{eqnarray}
where $R_0:=0$ and $T_{N+1}:=0$. We will now prove the following lemma.
%l13
\begin{lemma}\label{lem0705}
Suppose that $\beta=0$ and $p>1/2$. Let $M \in (0,\infty)$. There then
exists $C \in (0,\infty)$ such that for all but finitely
many $S \in \D$,
\[ \Exp_{0,p} [ f_2 (\xi_{t+1})^M - f_2 (\xi_{t})^M
 |  \xi_t = S ] \leq - C f_2(S)^{M-(5/6)}.\]
\end{lemma}

\begin{pf}
In this proof,\vspace*{1pt} we write
$\Delta_2(S):= \Exp_{0,p} [ f_2 (\xi_{t+1})^M - f_2 (\xi_{t})^M
 |  \xi_t = S ]$, which we calculate by summing over all the possible
exclusion moves. The $e_j^{10 \mapsto 01}$ transition has probability
$(1-p)/(2N+1)$ and changes $f_2(S)^M$ by
\[
\bigl(f_2(S)+1+R_j+T_{j+1} \bigr)^M - f_2(S)^M
= f_2(S)^M  \biggl[  \biggl( 1+ \frac{1+R_j + T_{j+1}}{f_2(S)}
\biggr)^M - 1  \biggr].
\]
Since $R_j + T_{j+1} \leq |S| = \mathrm{O}(f_2(S)^{1/2})$ for any $j$,
by (\ref{ineq1}),   Taylor's theorem  yields
\[  \bigl(f_2(S)+1+R_j+T_{j+1} \bigr)^M - f_2(S)^M
= f_2(S)^M  \biggl[ M \frac{1+R_j+T_{j+1}}{f_2(S)} + \mathrm{O} ( f_2 (S)^{-1} )
\biggr]. \]
Proceeding similarly for the $e_j^{01 \mapsto 10}$ transitions and
summing, we obtain
\begin{eqnarray}\label{0705d}
\Delta_2(S) &=& \frac{M}{2N+1} f_2(S)^{M-1}  \Biggl[
N+(1-p) + (1-2p) \sum_{j=1}^N (R_j +T_j)   \Biggr]\nonumber
\\[-8pt]\\[-8pt]
&&{}+\mathrm{O} (f_2(S)^{M-1}),\nonumber
\end{eqnarray}
where the implicit constant in $\mathrm{O}(\cdot)$ does not depend on $S$.
From (\ref{0705d}) with
(\ref{0825a}) and (\ref{ineq1}), we then obtain that for some $C_1,C_2 \in (0,\infty)$,
\[  \Delta_2(S) \leq   C_1 f_2(S)^{M-1} - C_2 f_2 (S)^{M-(5/6)}.\]
This yields the result.
\end{pf}

\begin{pf*}{Proof of Theorem~\ref{thmex1}}
Take $X_t = f_2 (\xi_t)^M$ for some $M\geq 1$.\vspace*{2pt} From Lemma~\ref{lem0705}, we then have that
$\Exp [ X_{t+1}- X_t  |  \xi_t = S ] \leq - C X_t^{1-{5}/{(6M)}}$
for all but finitely many $S$. We can\vspace*{1pt} thus apply Lemma~\ref{aimthm}
to obtain $\Exp [ \tau^{6M/5} ] < \infty$.
Since $M \geq 1$ was arbitrary, the theorem follows.
\end{pf*}

%s8 ###
\section{Size of the hybrid zone}\label{size}

We now prove the almost sure bounds on the rate
of growth of $| \xi_t |$ stated in Section~\ref{results}.
%l14
\begin{lemma} \label{lem4}
Let $\beta \in [0,1]$, $p \in [0,1]$. For any $\eps>0$,
$\Pr_{\beta,p}$-a.s., for all but finitely many $t$,
\[ \max_{0 \leq s \leq t} f_1 (\xi_s) \leq t (\log t)^{1+\eps}.\]
\end{lemma}
\begin{pf} From (\ref{0704c}),  we have that for any $S \in \D$,
\[ \Exp_{\beta,p} [ f_1 (\xi_{t+1} ) - f_1 (\xi_t)  |  \xi_t = S ] \leq
\frac{ N(S) +1}{2N(S) +1} \leq 1.\]
We can then apply Lemma~\ref{lem13}(ii)
with $X_t = f_1 (\xi_t)$ to obtain the result.
\end{pf}

\begin{pf*}{Proof of Theorem~\ref{thm3}}
Lemma~\ref{lem4}, with the simple inequality $f_1(\xi_t) \geq N(\xi_t)^2/2$,
implies the $p<1/2$ case of (\ref{yy0}).
By (\ref{0704a}), we have that for $p \geq 1/2$  and all $S \in \D$,
\[
\Exp_{\beta,p} [ f_2 (\xi_{t+1}) - f_2(\xi_t)  |  \xi_t = S ] \leq\frac{1-\beta}{2}.
\]
Hence, Lemma~\ref{lem13}(ii) with $X_t = f_2(\xi_t)$  yields,
for any $\eps>0$,
$\Pr_{\beta,p}$-a.s.,
\begin{eqnarray} \label{xx1} \max_{0 \leq s \leq t} f_2
(\xi_s)  \leq t (\log t)^{1+\eps}
\end{eqnarray}
for all but finitely many
$t$.  (\ref{yy1}) then follows from (\ref{xx1}) with (\ref{ineq1})
and the $p \geq 1/2$ case of (\ref{yy0}) follows  from (\ref{xx1})
with the simple inequality $f_2 (\xi_t) \geq (N(\xi_t))^3/3$
(obtained by replacing each $m_i$ and $n_i$ by $1$ in the
definition of $f_2$).
\end{pf*}

For the remainder of this section, we concentrate on the pure exclusion process,
that is,~when $\beta =0$. Again, the Lyapunov function $f_1$ will be
a primary tool here; the next result describes its behaviour in
this case. We use the abbreviation $N_t:= N(\xi_t)$.
%l15
\begin{lemma}
Let $\beta =0$, $p \in [0,1]$. Then $f_1(\xi_t)$ has
transition probabilities $p_j = \Pr_{0,p} ( f_1 (\xi_{t+1}) -
f_1(\xi_t) = j  |  \F_t )$ for jumps $j \in \{-1,0,+1\}$,
where $p_{-1}+p_0+p_1=1$ and
%e9 ###
\begin{equation}\label{eq33.5}
p_{-1} = p \frac{N_t}{2N_t+1} \leq \frac{p}{2},\qquad
p_0 = \frac{N_t+p}{2N_t+1},\qquad
p_1 = (1-p)  \frac{N_t+1}{2N_t+1} \geq \frac{1-p}{2}.
\end{equation}
Hence, for all $t \in \Z^+$,
\begin{eqnarray} \label{zz4}
f_1 (\xi_t) \leq
f_1 (\xi_0) + t.
\end{eqnarray} Moreover, when $p <1/2$,
for any $c \in (0, (1/2)-p)$,  $\Pr_{0,p}$-a.s.,
for all but finitely many $t$,
\begin{eqnarray} \label{zz3} f_1 (\xi_t) \geq ct.
\end{eqnarray}
\end{lemma}
\begin{pf} (\ref{eq33.5}) follows from equations~(5.5)  and (5.6) in
\cite{bfmp}.  (\ref{zz4}) is then immediate.
From (\ref{eq33.5}), we have that $\xi_t$
stochastically dominates $\xi_0+\sum_{s=1}^t W_s$,
where $W_1, W_2,\ldots$ are i.i.d.~random
variables taking values $+1$, $0$, $-1$
with probabilities $q/2$, $1/2$, $p/2$, respectively.
Hence, the SLLN and the fact that $\Exp [ W_1]
=(1/2)-p$ yields (\ref{zz3}) for $p<1/2$.
\end{pf}
%cor2
\begin{corollary} Suppose that $\beta=0$ and $p < 1/2$. There then
exists $c'(p)>0$ such that
for any $c \in (0,c'(p))$,
$\Pr_{0,p}$-a.s., for all but finitely many $t$,
\begin{eqnarray}\label{zz1}
|\xi_t| \geq c t^{1/2}.\end{eqnarray}
Suppose that $\beta =0$.
There then exists $C \in (0,\infty)$ such that for
any $p \in [0,1]$, $\Pr_{0,p}$-a.s.,
for all but finitely many $t \in \Z^+$,
\begin{eqnarray}\label{zz2} N(\xi_t) \leq C t^{1/2}.\end{eqnarray}
\end{corollary}
\begin{pf}
The bound (\ref{zz1}) follows from (\ref{zz3}) together with
(\ref{ineq1}); (\ref{zz2}) follows from (\ref{zz4}) with the simple
inequality $f_1 (\xi_t) \geq (N(\xi_t))^2/2$.
\end{pf}

The next two lemmas give some properties of the process
$(|\xi_t|)_{t \in \Z^+}$. Recall the definition of configuration
$\D_1$ from (\ref{d1def}).
%l16
\begin{lemma}  \label{l1}
Suppose that $\beta=0$ and $p \in [0,1]$.   For any $t \in\Z^+$,
we have that
\begin{eqnarray}    \label{tt1}
\Pr_{0,p} ( |\xi_{t+1}| = 2  |  \xi_t = \D_0 )
&= &1 - \Pr_{0,p} ( |\xi_{t+1}| = 0  |  \xi_t = \D_0 ) = 1-p\quad\mbox{and} \\
\label{tt1a}
  \Pr_{0,p} ( |\xi_{t+1}| = j  |  \xi_t = \D_1 ) &=& \frac{p}{3},
\frac{1+p}{3}, \frac{2(1-p)}{3}  \qquad\mbox{for } j = 0, 2, 3, \mbox{ respectively.}
\end{eqnarray}
For any $t \in\Z^+$, conditional on $\xi_t \in \D \setminus \{ \D_0, \D_1 \}$,
$|\xi_{t+1}|-|\xi_t|$ takes values only in $\{-1,0,+1\}$ and
for any $S \in \D \setminus \{ \D_0, \D_1 \},$
\begin{eqnarray}    \label{tt2}
\Pr_{0,p} ( |\xi_{t+1}| - |\xi_t| = 1  |  \xi_t = S )
& =& \frac{2(1-p)}{2N(S)+1} \quad\mbox{and}\\
\label{tt3}
\Pr_{0,p} ( |\xi_{t+1} | - |\xi_t| = -1  |  \xi_t = S )
& =& \frac{p( \1_{\{n_1(S) =1\}} + \1_{\{m_{N(S)}(S) =1\}} )}{2N(S)+1}
\leq \frac{2p}{2N(S)+1}.
\end{eqnarray}
\end{lemma}
\begin{pf}
The statements (\ref{tt1}) and (\ref{tt1a}) are straightforward. Suppose that
$\xi_t = S$ for some $S \in \D \setminus \{ \D_0, \D_1\}$. Then $|S| \geq 2$
and exclusion moves cannot effect a change of magnitude more than~$1$.
We have that $|\xi_{t+1}| = |S| +1$
if and only if we select (with probability $2/(2N(S)+1)$)
one of the two extreme $10$ pairs and then
(with probability $1-p$) we flip the $10$ to a $01$. Similarly,
$|\xi_t|$ can decrease by $1$ if and only if there exists
a configuration $\ldots 11101\ldots$ at the left end
or a configuration $\ldots 01000\ldots$ at the right end, and then
we select the $01$ and flip to $10$.
The statement of the lemma follows.
\end{pf}

%l17
\begin{lemma}
If $\beta=0$ and $p=1/2$, then
\begin{eqnarray} \label{vv2}
\Exp_{0,1/2} [ |\xi_{t+1}|^3 - |\xi_t|^3  |  \mathcal{F}_t ] \geq 4 \qquad\mbox{a.s.}
\end{eqnarray}
\end{lemma}
\begin{pf} First, we have from (\ref{tt1}) and (\ref{tt1a})
that
\[
\Exp_{0,1/2} [ |\xi_{t+1}|^3 - |\xi_t|^3  |  \xi_t = \D_0 ]= 4,\qquad
\Exp_{0,1/2} [ |\xi_{t+1}|^3 - |\xi_t|^3  |  \xi_t = \D_1 ] = 5.
\]
It thus remains to consider the case where $\xi_t =S$ for $S
\in \D \setminus \{ \D_0,\D_1 \}$.
Here, from (\ref{tt2}) and (\ref{tt3}), writing $N = N(S)$, we have
\begin{eqnarray*}
\Exp_{0,1/2} [ |\xi_{t+1}|^3 - |\xi_t|^3  |  \xi_t = S ]
&\geq&
\frac{1}{2N+1}  \bigl[  \bigl( (|S|+1)^3 - |S|^3  \bigr) +  \bigl(
(|S|-1)^3 - |S|^3  \bigr)  \bigr] \\
&= &\frac{6 |S|}{2N+1}
\geq \frac{12N}{2N+1}   \geq 4
\end{eqnarray*}
since $|S| \geq 2N(S)$ and $N(S) \geq 1$ for all $S \neq \D_0$.
\end{pf}

\begin{pf*}{Proof of Theorem~\ref{thm5}}
The upper bound in
the theorem is implied by (\ref{yy1}). For the lower bound, use
Theorem 3.3 of \cite{mvw} with, in the notation of that paper,
$f(x) = x^3$ and $Y_n = |\xi_n|$.
Using (\ref{vv2}) and the fact that
$|\xi_t|$ has uniformly bounded jumps (see Lemma~\ref{l1}), we then
obtain
the desired result.
\end{pf*}

We now work toward the upper bound for $|\xi_t|$, for $p \in [0,1]$,
given in Theorem~\ref{twothirds}.
Define the function $\rho^2$ by
$\rho^2 (\D_0):= 0$ and, for $S \in \D \setminus \{ \D_0\},$
\begin{eqnarray}\label{rhodef}
\rho^2(S):= \sum_{i=1}^N m_i^2 + \sum_{i=1}^N n_i ^2.
\end{eqnarray}
%l18
\begin{lemma}
For any $S \in \D \setminus\{ \D_0 \}$,
\begin{eqnarray}\label{ww1c}
|S| \leq \frac{1}{2N(S)} |S|^2 \leq \rho^2 (S) \leq |S|^2.
\end{eqnarray}
\end{lemma}
\begin{pf}
Suppose that $S \in \D \setminus \{ \D_0\}$.  For the upper bound, we have
\(  \rho^2 (S) \leq \sum_{i=1}^N (m_i +n_i)^2  \leq |S|^2 \).
For the lower bound, \( \frac{1}{N} \sum_{i=1}^N m_i^2 \geq  (
\frac{1}{N} \sum_{i=1}^N m_i  )^2 \), from Jensen's inequality,
and similarly for the $n_i$. Hence,
\[
\rho^2 (S) \geq \frac{1}{N} ( R_N^2 + T_1^2) \geq \frac{1}{2N} (R_N
+ T_1 )^2 = \frac{1}{2N} |S|^2 \geq |S|
\]
since $|S| \geq 2N$, completing the proof.
\end{pf}
%l19
\begin{lemma}
Suppose that $\beta=0$ and $p \in [0,1]$. For $t \in \Z^+$, we have
\begin{eqnarray} \label{ccc1}
\Exp_{0,p} [ \rho^2(\xi_{t+1}) -\rho^2 (\xi_t)  |  \F_t ] \leq 2 \qquad\mbox{a.s.};
\end{eqnarray}
moreover, for any $\eps>0$, $\Pr_{0,p}$-a.s., for all
but finitely many $t \in \Z^+$,
\begin{eqnarray} \label{ccc2}
\max_{0 \leq s \leq t} \rho^2 (\xi_s) \leq t (\log t)^{1+\eps}.
\end{eqnarray}
\end{lemma}
\begin{pf}
We start by proving (\ref{ccc1}).
First, we note that for any $t \in \Z^+$,
\[ \Exp_{0,p} [ \rho^2 (\xi_{t+1}) - \rho^2 (\xi_t)  |  \xi_t = \D_0 ] = 2(1-p)
\leq 2. \]
We next need to verify (\ref{ccc1}) for any
configuration $S \in \D \setminus \{ \D_0\}$.

Let $\Delta_{1,i} (S)$ denote the change in $\rho^2 (S)$ when a
$01\mapsto 10$ exclusion move is performed on the $i$th $01$ pair in
$S$ ($i=1,\ldots,N$). Similarly, let $\Delta_{2,i} (S)$ denote the
change in $\rho^2 (S)$ when a $10\mapsto 01$ exclusion move is
performed on the $i$th $10$ pair ($i=0,\ldots,N$). Thus,
\[ \Delta_{1,i} (S):= \rho^2 ( e_i^{01 \mapsto 10} (S))
- \rho^2 (S),\qquad
\Delta_{2,i} (S):= \rho^2 ( e_i^{10 \mapsto 01} (S))
- \rho^2 (S),\]
in the notation of
Section~\ref{config}. We then have
\begin{eqnarray}\label{kk1}
\Exp_{0,p} [ \rho^2 (\xi_{t+1}) - \rho^2 (\xi_t)  |
\xi_t = S ] = \frac{1}{2N+1}  \Biggl( p \sum_{i=1}^N
\Delta_{1,i} (S) + q \sum_{i=0}^{N}
\Delta_{2,i} (S)  \Biggr).
\end{eqnarray}
We compute the two sums on the right-hand side
of (\ref{kk1}) separately. First, consider all $N$ possible
exclusion moves $01\mapsto 10$. Separating out the cases when $m_i=1$
or $n_i=1$,
\begin{eqnarray*}
\Delta_{1,i}(S)  =  -(2m_i-2)\1_{\{m_i>1\}}-(2n_i-2)\1_{\{n_i>1\}}
+2m_{i-1} \1_{\{n_i=1\}}+2n_{i+1} \1_{\{m_i=1\}},
\end{eqnarray*}
with the convention that $n_{N+1}=m_{0}=-1/2$ to make this
formula correct for $i=1$ and $i=N$. Since $(x-1) \1_{\{ x > 1\}} = x-1$
for $x \in \N$,  this last equation is
\begin{eqnarray*}
\Delta_{1,i}(S) = 2\bigl[ 2-m_i-n_i+m_{i-1} \1_{\{n_i=1\}}+n_{i+1} \1_{\{m_i=1\}}\bigr].
\end{eqnarray*}
Hence, summing over $i \in \{1,\ldots,N\}$ gives
\begin{eqnarray}  \label{kk2}
\hspace*{-12pt}\frac{1}{2} \sum_{i=1}^N \Delta_{1,i} (S)
= 2N -\sum_{i=1}^{N}  (m_i+n_i)+\sum_{i=0}^{N-1} m_{i}
\1_{\{n_{i+1}=1\}} + \sum_{i=2}^{N+1} n_{i} \1_{\{m_{i-1}=1\}}
\leq 2N.
\end{eqnarray}
Similarly, a  $10\mapsto 01$ exclusion move on the $i$th $10$ pair
($i=0,1,\ldots,N$) contributes
\[ \Delta_{2,i}(S) =
2 \bigl[ 2 -m_i-n_{i+1} +m_{i+1} \1_{\{n_{i+1}=1\}}+n_{i}
\1_{\{m_i=1\}} \bigr],
\]
with the conventions that $n_0=m_{N+1}=0$ and $n_{N+1}=m_0=1/2$ to make this
formula correct for $i=0$ and $i=N$. Summing, as before,
\begin{eqnarray}
\label{kk3}
\frac{1}{2}
\sum_{i=0}^{N} \Delta_{2,i} (S) = 2N +1  -\sum_{i=1}^{N} m_i \1_{\{n_{i}>1\}}
-\sum_{i=1}^{N} n_{i} \1_{\{m_{i}>1\}} \leq 2N+1.
\end{eqnarray}
Combining (\ref{kk2}) and (\ref{kk3}) with (\ref{kk1}), we conclude that
\begin{eqnarray*}
\Exp_{0,p} [ \rho^2(\xi_{t+1}) -\rho^2 (\xi_t)  |  \xi_t = S ]
\leq \frac {2 (2N +q)}{2N+1} \leq 2,
\end{eqnarray*}
which is (\ref{ccc1}). Finally, (\ref{ccc2}) follows from (\ref{ccc1})
with Lemma~\ref{lem13}(ii), taking  $X_t = \rho^2 (\xi_t)$.
\end{pf}

Suppose that $p \in [0,1]$. Then, from (\ref{ccc2}), (\ref{zz2}) and the
middle inequality in (\ref{ww1c}), we obtain
an upper bound for $\max_{0 \leq s \leq t} |\xi_s|$ of order $t^{3/4}$
(ignoring logarithmic factors).
In order to prove the upper bound in
Theorem~\ref{twothirds}(ii), we will
give an argument that improves the $3/4$ to $2/3$.
We start with a simple inequality.
%l20
\begin{lemma}\label{stas}
Let $N \in \N$. Suppose that $n_1,n_2,\ldots, n_N \geq 0$. If, for
some $A, B >0$,
\[  \sum_{i=1}^N n_i^2\le A \quad\mbox{and}\quad \sum_{i=1}^N i n_i
\le B,\quad \mbox{then}\quad
\sum_{i=1}^N n_i \le (6AB)^{1/3}.
\]
\end{lemma}
\begin{pf}
In the elementary inequality
\(   (\sum_{i=1}^N  n_i  )^3 \leq  3\sum_{i=1}^N  n_i^2 \sum_{i=1}^N
n_i + 3 \sum_{i=1}^N n_i  (\sum_{j=1}^{i-1} n_j  )^2 \),
apply Jensen's inequality to the final term to obtain
\[
\Biggl(\sum_{i=1}^N  n_i  \Biggr)^3 \leq
3 \sum_{i=1}^N  n_i^2 \sum_{i=1}^N i
n_i  + 3 \sum_{i=1}^N i n_i  \sum_{j=1}^{i-1} n_j^2
\le 6AB.
\]
\upqed
\end{pf}

\begin{pf*}{Proof of Theorem~\ref{twothirds}}
Part (i) of the theorem is (\ref{zz2}) and the lower bound in part
(ii) of the theorem is (\ref{zz1}). We now derive the upper bound in part
(ii). Since each block of $1$'s has at least one element,
observe that $\sum_{i=1}^N n_i (N-i) \le  f_1(S)$
and also that $\sum_{i=1}^N n_i^2  \le  \rho^2(S)$.
Lemma~\ref{stas} thus implies that for any $S \in \D \setminus \{\D_0\}$,
\[ \sum_{i=1}^{N} n_i \leq (6 f_1 (S) \rho^2 (S) )^{1/3} \leq 2 ( f_1
(S) \rho^2 (S) )^{1/3},\]
and the same argument applies for $\sum_{i=1}^N m_i$. Hence,
for any $S \in \D$,
\begin{eqnarray}\label{aaa}
|S| \leq 4 (f_1 (S) \rho^2 (S) )^{1/3}.
\end{eqnarray}
Taking $S=\xi_t$, we have
$f_1(\xi_t)  \le C_1 t$ for all $t$
and some $C_1 \in (0,\infty)$,
by
(\ref{zz4}). Also, for any $\eps>0$,
$\Pr_{0,p}$-a.s., $\rho^2(\xi_t) \le C_2 t
(\log t)^{1+\eps}$ for all $t$, by (\ref{ccc2}), for some $C_2 \in (0,\infty)$.
Using these bounds in (\ref{aaa}) completes the proof.
\end{pf*}

\section*{Acknowledgements}

This work was done while Andrew Wade was
at the University of Bristol,
supported by the Heilbronn Institute for Mathematical Research.
The authors are grateful to an anonymous referee
for a careful reading of the paper and for
drawing our attention to some relevant references.

\printhistory

\end{document}